\input amstex
\documentstyle{amsppt}
\magnification=\magstep 1
\newcount\zero \zero=0
\newcount\prelim \prelim=1
\newcount\comF \comF=2
\newcount\comM \comM=3
\newcount\xact \xact=4
\newcount\misc \misc=5
\define\Adj{\operatorname{Adj}}
\predefine\Hunumlaut{\H}
\redefine\H{\operatorname{H}}
\define\a{\alpha}
\predefine\barunder{\b}
\redefine\b{\beta}
\predefine\cedilla{\c}
\define\g{\gamma}
\predefine\dotunder{\d}
\redefine\d{\delta}
\define\e{\varepsilon}
\define\f{\varphi}
\define\({\left(}
\define\){\right)}
\define\[{\left[}
\define\]{\right]}
\define\Tor{\operatorname{Tor}}
\define\Ext{\operatorname{Ext}}
\define\w{\wedge}
\define\W{\tsize\bigwedge}
\predefine\tieaccent{\t}
\redefine\t{\otimes}
\define\iso{\cong}
\define\id{\text{id}}
\define\grade{\operatorname{grade}}
\define\pd{\operatorname{pd}}
\define\depth{\operatorname{depth}}
\define\rank{\operatorname{rank}}
\define\Ker{\operatorname{Ker}}
\define\Hom{\operatorname{Hom}}
\define\Z{\Bbb Z}
\define\F{\Bbb F}
\define\G{\Bbb G}
\define\C{\Bbb C}
\define\im{\operatorname{im}}
\define\p{\oplus}
\define \maxm {\frak m}
\def\M{\Bbb M}
\define\fakeht{\vphantom{E^{E^{E}}_{E_{E}}}}
\def\maxM{\frak M}
\def\diag{\operatorname{diag}}
\def\u{\bold u}
\def\v{\bold v}
\def\X{\bold X}
\def\bx{\sigma}
\predefine\image{\Im}
\redefine\Im{\operatorname{Im}}
\def\N{\Bbb N}

\def\Nh{\widehat{\N}}
\predefine\ELL{\L}
\redefine\L{\Bbb L}
\def\incl{\operatorname{incl}}
\predefine\ETA{\eta}
\redefine\eta{e_{n}}

\def\Ft{\widetilde{\F}}
\def\ft{\widetilde{f}}
\predefine\BBB{\B}
\redefine\B{\Bbb B}
\predefine\DDD{\D}
\redefine\D{\Bbb D}
\predefine\TTT{\T}
\redefine\T{\Bbb T}
\define\proj{\operatorname{proj}}
\def\A{\Bbb A}
\def\nat{\operatorname{nat}}
\def\Ah{\widehat{\A}}
\def\Bh{\widehat{\B}}
\def\Ch{\widehat{\C}}
\def\Dh{\widehat{\D}}
\def\Fh{\widehat{\F}}
\def\fh{\widehat{f}}
\def\tri{\frak f}
\define\triv{\ltimes}
\define\Fakeht{\vphantom{\matrix X\\X\endmatrix}}

\chardef\oldatsign=\catcode`\@
\catcode`\@=11
\newif\ifdraftmode			% New- print marginal notes if true
\global\draftmodefalse

% End of \draft.


%===============================================================================
% Font loading. Normal text size is 12pt.
%
\font@\twelverm=cmr12 % roman text 
\font@\twelvei=cmmi12 \skewchar\twelvei='177 % math italic
\font@\twelvesy=cmsy10 scaled\magstep1 \skewchar\twelvesy='060 % math symbols
\font@\twelveex=cmex10 scaled\magstep1 % math extensions
\font@\twelvemsa=msam10 scaled\magstep1 % AMS extra symbols
\font@\twelvemsb=msbm10 scaled\magstep1 % AMS extra symbols
\font@\twelvebf=cmbx12 % boldface extended
\font@\twelvett=cmtt12 % typewriter
\font@\twelvesl=cmsl12 % slanted roman
\font@\twelveit=cmti12 % text italic
\font@\twelvesmc=cmcsc10 scaled\magstep1 % caps and small caps
%
%     the normal script font is nine point
%
\font@\ninerm=cmr9 % roman text
\font@\ninei=cmmi9 \skewchar\ninei='177 % math italic
\font@\ninesy=cmsy9 \skewchar\ninesy='60 % math symbols
\font@\ninemsa=msam9
\font@\ninemsb=msbm9
\font@\ninebf=cmbx9
%
%	the normal scriptscript font is 7 point; already loaded by amstex
%	or plain.
%
%     some specialty fonts
%
\font@\ttlrm=cmbx12 scaled \magstep2 % a title font
\font@\ttlsy=cmsy10 scaled \magstep3 % for \AmSTeX in titles
\font@\tensmc=cmcsc10 % caps and small caps
 % a little caps type
 % a little tt type
%
%     twelve point is used for normal printing.
%
\def\normaltype{%   twelve point stuff
	\def\pointsize@{12}%
	\abovedisplayskip18\p@ plus5\p@ minus9\p@
	\belowdisplayskip18\p@ plus5\p@ minus9\p@
	\abovedisplayshortskip1\p@ plus3\p@
	\belowdisplayshortskip9\p@ plus3\p@ minus4\p@
	\textonlyfont@\rm\twelverm
	\textonlyfont@\it\twelveit
	\textonlyfont@\sl\twelvesl
	\textonlyfont@\bf\twelvebf
	\textonlyfont@\smc\twelvesmc
	\ifsyntax@
		\def\big##1{{\hbox{$\left##1\right.$}}}%
	\else
		\let\big\twelvebig@
 \textfont0=\twelverm \scriptfont0=\ninerm \scriptscriptfont0=\sevenrm
 \textfont1=\twelvei  \scriptfont1=\ninei  \scriptscriptfont1=\seveni
 \textfont2=\twelvesy \scriptfont2=\ninesy \scriptscriptfont2=\sevensy
 \textfont3=\twelveex \scriptfont3=\twelveex  \scriptscriptfont3=\twelveex
 \textfont\itfam=\twelveit \def\it{\fam\itfam\twelveit}%
 \textfont\slfam=\twelvesl \def\sl{\fam\slfam\twelvesl}%
 \textfont\bffam=\twelvebf \def\bf{\fam\bffam\twelvebf}%
 \scriptfont\bffam=\ninebf \scriptscriptfont\bffam=\sevenbf
 \textfont\ttfam=\twelvett \def\tt{\fam\ttfam\twelvett}%
 \textfont\msafam=\twelvemsa \scriptfont\msafam=\ninemsa
 \scriptscriptfont\msafam=\sevenmsa
 \textfont\msbfam=\twelvemsb \scriptfont\msbfam=\ninemsb
 \scriptscriptfont\msbfam=\sevenmsb
	\fi
 \normalbaselineskip=\twelvebaselineskip
 \setbox\strutbox=\hbox{\vrule height12\p@ depth6\p@
      width0\p@}%
 \normalbaselines\rm \ex@=.2326ex%
}% End of \normaltype.
%
%
%     ten point is used for smalltype
%
\def\smalltype{%   ten point stuff
	\def\pointsize@{10}%
	\abovedisplayskip12\p@ plus3\p@ minus9\p@
	\belowdisplayskip12\p@ plus3\p@ minus9\p@
	\abovedisplayshortskip\z@ plus3\p@
	\belowdisplayshortskip7\p@ plus3\p@ minus4\p@
	\textonlyfont@\rm\tenrm
	\textonlyfont@\it\tenit
	\textonlyfont@\sl\tensl
	\textonlyfont@\bf\tenbf
	\textonlyfont@\smc\tensmc
	\ifsyntax@
		\def\big##1{{\hbox{$\left##1\right.$}}}%
	\else
		\let\big\tenbig@
	\textfont0=\tenrm \scriptfont0=\sevenrm \scriptscriptfont0=\fiverm 
	\textfont1=\teni  \scriptfont1=\seveni  \scriptscriptfont1=\fivei
	\textfont2=\tensy \scriptfont2=\sevensy \scriptscriptfont2=\fivesy 
	\textfont3=\tenex \scriptfont3=\tenex \scriptscriptfont3=\tenex
	\textfont\itfam=\tenit \def\it{\fam\itfam\tenit}%
	\textfont\slfam=\tensl \def\sl{\fam\slfam\tensl}%
	\textfont\bffam=\tenbf \def\bf{\fam\bffam\tenbf}%
	\scriptfont\bffam=\sevenbf \scriptscriptfont\bffam=\fivebf
	\textfont\msafam=\tenmsa
	\scriptfont\msafam=\sevenmsa
	\scriptscriptfont\msafam=\fivemsa
	\textfont\msbfam=\tenmsb
	\scriptfont\msbfam=\sevenmsb
	\scriptscriptfont\msbfam=\fivemsb
		\textfont\ttfam=\tentt \def\tt{\fam\ttfam\tentt}%
	\fi
 \normalbaselineskip 14\p@
 \setbox\strutbox=\hbox{\vrule height10\p@ depth4\p@ width0\p@}%
 \normalbaselines\rm \ex@=.2326ex%
}% End of \smalltype.

\def\titletype{%   title fonts
	\def\pointsize@{17}%
	\textonlyfont@\rm\ttlrm
	\ifsyntax@
		\def\big##1{{\hbox{$\left##1\right.$}}}%
	\else
		\let\big\twelvebig@
		\textfont0=\ttlrm \scriptfont0=\twelverm
		\scriptscriptfont0=\tenrm
		\textfont2=\ttlsy \scriptfont2=\twelvesy
		\scriptscriptfont2=\tensy
	\fi
	\normalbaselineskip 25\p@
	\setbox\strutbox=\hbox{\vrule height17\p@ depth8\p@ width0\p@}%
	\normalbaselines
	\rm
	\ex@=.2326ex%
}% End of \titletype.

\def\tenbig@#1{% Same as PLAIN definition of \big .
	{%
		\hbox{%
			$%
			\left
			#1%
			\vbox to8.5\p@{}%
			\right.%
			\n@space
			$%
		}%
	}%
}% End of \tenbig@.

\def\twelvebig@#1{%
	{%
		\hbox{%
			$%
			\left
			#1%
			\vbox to10.2\p@{}%  This is just a guess
			\right.%
			\n@space
			$%
		}%
	}%
}% End of \twelvebig@.

%==============================================================================
%	Macros for symbolic numbering of theorems, figures, etc.
%
%	Purpose:	To provide multiple tracks of symbolic labels for
%			theorems, figures, etc., which will automatically be
%			converted into numbers in a specified format.
%
%	A 'track' or 'type' of labels is identified by a string of letters,
%	with case distinguished.  A symbolic label is a string of any
%	characters.
%
%	Forward references are resolved using two passes and an external file,
%	\jobname.xref.
%
%	Associated with each symbolic label, there are two macros, one holding
%	the actual label, and the other holding the 'state' of the label, which
%	is used for error-checking.
%
%==============================================================================
\newif\ifl@beloutopen
\newwrite\l@belout
\newread\l@belin

\global\let\currentfile=\jobname

% getfile
% =======
% #1: tex file to input
\def\getfile#1{%
	\immediate\closeout\l@belout
	\global\l@beloutopenfalse
	\gdef\currentfile{#1}%
	\input #1%
	\par
	\newpage
}% End of \getfile.

% getxrefs
% ========
% #1: list of other files, separated by commas, from which to grab labels.
\def\getxrefs#1{%
	\bgroup
		\def\gobble##1{}% used to discard \end when we're done
		\edef\list@{#1,}%
		\def\gr@boff##1,##2\end{% process 1st thing in list
			\openin\l@belin=##1.xref
			\ifeof\l@belin
			\else
				\closein\l@belin
				\input ##1.xref
			\fi
			\def\list@{##2}%
			\ifx\list@\empty
				\let\next=\gobble
			\else
				\let\next=\gr@boff
			\fi
			\expandafter\next\list@\end
		}%
		\expandafter\gr@boff\list@\end
	\egroup
}% End of \getxrefs.

% testdefined
% ===========
% #1: control sequence to test
% #2: stuff to execute if #1 is defined
% #3: stuff to execute if #1 is not defined
\def\testdefined#1#2#3{%
	\expandafter\ifx
	\csname #1\endcsname
	\relax
	#3%
	\else #2\fi
}% End of \testdefined

\def\document{%
%	\openout\contents=\jobname.contents
%	\normaltype
	\minaw@11.11128\ex@ % minimum arrow with for @>...>...> macro
%	\pageno=-2 % for preliminary pages; -1 will be for table of cont.
	\def\alloclist@{\empty}%
	\def\fontlist@{\empty}%
	\openin\l@belin=\jobname.xref	% Input \jobname.xref if it exists.
	\ifeof\l@belin\else
		\closein\l@belin
		\input \jobname.xref
	\fi
}% End of \document.

% getst@te
% ========
% #1: type of label
% #2: symbolic label
\def\getst@te#1#2{%
	\edef\st@te{\csname #1s!#2\endcsname}%
	\expandafter\ifx\st@te\relax
		\def\st@te{0}%
	\fi
}% End of \getst@te

% setst@te
% ========
% #1: type of label
% #2: symbolic label
% #3: new state
\def\setst@te#1#2#3{%
	\expandafter
	\gdef\csname #1s!#2\endcsname{#3}%
}% End of setst@te.

% setupautolabel
% ==============
% #1: type of label
% #2: form of actual label (a token list, which should involve the count
%	register \#1Number).
\outer\def\setupautolabel#1#2{%
	\def\newcount@{\global\alloc@0\count\countdef\insc@unt}	% Normally,
		% \newcount is \outer, so we have to repeat the definition
		% here, using some macros from plain.
	\def\newtoks@{\global\alloc@5\toks\toksdef\@cclvi}% See above.
	\expandafter\newcount@\csname #1Number\endcsname
	\expandafter\global\csname #1Number\endcsname=1%
	\expandafter\newtoks@\csname #1l@bel\endcsname
	\expandafter\global\csname #1l@bel\endcsname={#2}%
}% End of \setupautolabel.

% reflabel
% ========
% #1: type of label
% #2: symbolic label
\def\reflabel#1#2{%
	\testdefined{#1l@bel}% See whether the type is known.
	{%	Yes, the type is known.
		\getst@te{#1}{#2}%
		\ifcase\st@te
			%	State 0: the label is undefined.
			???%	Insert dummy label.
			\message{Unresolved forward reference to
				label #2. Use another pass.}%
		\or	%	state 1: defined by file, not yet referenced
			\setst@te{#1}{#2}2%       Why do I need this %?
			\csname #1l!#2\endcsname % Insert the label.
		\or	%	State 2: defined by file, referenced
			\csname #1l!#2\endcsname % Insert the label.
		\or	%	State 3: defined by a \setlabel
			\csname #1l!#2\endcsname % Insert the label.
		\fi
	}{%	No, the type is unknown.
		{\escapechar=-1 % to make the next statement work
		\errmessage{You haven't done a
			\string\\setupautolabel\space for type #1!}%
		}%
	}%
}% End of \reflabel.

{\catcode`\{=12 \catcode`\}=12
	\catcode`\[=1 \catcode`\]=2
	\xdef\Lbrace[{]%	Needed in order to write braces to a file.
	\xdef\Rbrace[}]%
]%

% setlabel
% ========
% #1: type of label
% #2: symbolic label
\def\setlabel#1#2{%
	\testdefined{#1l@bel}%	See whether the type is known.
	{%	Yes, the type is known.
		\edef\templ@bel@{\expandafter\the
			\csname #1l@bel\endcsname}%
		\def\@rgtwo{#2}%
		\ifx\@rgtwo\empty
		\else
			\ifl@beloutopen\else
				\immediate\openout\l@belout=\currentfile.xref
				\global\l@beloutopentrue
			\fi
			\getst@te{#1}{#2}%
			\ifcase\st@te
				% state 0: undefined
			\or	% state 1: defined by file, not yet used
			\or	% state 2: defined by file, referenced
				\edef\oldnumber@{\csname #1l!#2\endcsname}%
				\edef\newnumber@{\templ@bel@}%
				\ifx\newnumber@\oldnumber@
				\else
					\message{A forward reference to label 
						#2 has been resolved
						incorrectly.  Use another
						pass.}%
				\fi
			\or	% state 3: defined by \setlabel
				\errmessage{Same label #2 used in two
					\string\setlabel s!}%
			\fi
			\expandafter\xdef\csname #1l!#2\endcsname
				{\templ@bel@}%	Set the value of the label.
			\setst@te{#1}{#2}3%
			\immediate\write\l@belout % Save label value.
				{\string\expandafter\string\gdef
				\string\csname\space #1l!#2%
				\string\endcsname
				\Lbrace\templ@bel@\Rbrace
				}%
			\immediate\write\l@belout % Save label state.
				{\string\expandafter\string\gdef
				\string\csname\space #1s!#2%
				\string\endcsname
				\Lbrace 1\Rbrace
				}%
		\fi
		\templ@bel@	% Insert the label value.
		\expandafter\ifx\envir@end\endref % inside \ref?
			\gdef\marginalhook@{\marginal{#2}}%
		\else
			\marginal{#2}%    write symbolic label in margin
		\fi
		\expandafter\global\expandafter\advance	% Increment the counter
			\csname #1Number\endcsname
			by 1 %
	}{%	No, the type is unknown.
		{\escapechar=-1
		\errmessage{You haven't done a \string\\setupautolabel\space
			for type #1!}%
		}%
	}%
}% End of \setlabel.

%====================End of symbolic labelling macros==========================
%	Macro for labelling theorems, definitions, etc.
% The following macros allow forward references, at the expense of requiring
% two passes. 

\newcount\SectionNumber
\setupautolabel{t}{\number\SectionNumber.\number\tNumber}
\setupautolabel{r}{\number\rNumber}
\setupautolabel{T}{\number\TNumber}

\define\rref{\reflabel{r}}
\define\tref{\reflabel{t}}

\define\tnum{\setlabel{t}}
\define\rnum{\setlabel{r}}

%%%%%%%%%%%%%%%%%%%%%%%%%%%%%%%%%%%%%%%%%%%%%%%%%%%%%%%%%%%%%%%%%%%%%%%
%  macros for marginal notes
%
\def\strutdepth{\dp\strutbox}%
\def\strutheight{\ht\strutbox}%

\newif\iftagmode
\tagmodefalse

\let\old@tagform@=\tagform@
\def\tagform@{\tagmodetrue\old@tagform@}

\def\marginal#1{%
	\ifvmode
	\else
		\strut
	\fi
	\ifdraftmode
		\ifmmode
			\ifinner
				\let\Vorvadjust=\Vadjust
			\else%				display math mode
				\let\Vorvadjust=\vadjust
			\fi
		\else
			\let\Vorvadjust=\Vadjust
		\fi
		\iftagmode	% special case - tag of an equation
			\llap{%
				\smalltype
				\vtop to 0pt{%
					\pretolerance=2000
					\tolerance=5000
					\raggedright
					\hsize=.72in
					\parindent=0pt
					\strut
					#1%
					\vss
				}%
				\kern.08in
				\iftagsleft@
				\else
					\kern\hsize
				\fi
			}%
		\else% not tagmode
			\Vorvadjust{%
				\kern-\strutdepth % back up to baseline
				{%
					\smalltype
					\kern-\strutheight % match next baseline
					\llap{%
						\vtop to 0pt{%
							\kern0pt
							\pretolerance=2000
							\tolerance=5000
							\raggedright
							\hsize=.5in
							\parindent=0pt
							\strut
							#1%
							\vss
						}%
						\kern.08in
					}%
					\kern\strutheight
				}%
				\kern\strutdepth
			}% end of Vorvadjust.
		\fi% end iftagmode.
	\fi
}% End of \marginal.

% When a \vadjust occurs inside nested boxes, it doesn't seem to do anything,
% so we need to do some trickery to get the adjustment outside.

\newbox\Vadjustbox

\def\Vadjust#1{% My generalization of \vadjust
	\global\setbox\Vadjustbox=\vbox{#1}%
	\ifmmode
		\ifinner
			\innerVadjust
		\fi		%	don't do it in display math mode
	\else
		\innerVadjust
	\fi
}% End of \Vadjust.

\def\innerVadjust{%
	\def\nexti{\aftergroup\innerVadjust}%
	\def\nextii{%
		\ifvmode
			\hrule height 0pt % to prevent \baselineskip glue
			\box\Vadjustbox
		\else
			\vadjust{\box\Vadjustbox}%
		\fi
	}%
	\ifinner
		\let\next=\nexti
	\else
		\let\next=\nextii
	\fi
	\next
}%

\global\let\marginalhook@\empty

\catcode`\@=\oldatsign
 
\NoBlackBoxes
\topmatter
\keywords Acyclicity lemma, Exterior algebra, Finite free resolution, Gorenstein ideal, Koszul complex,  Multilinear algebra, Perfect ideal, Tor-algebra \endkeywords
\abstract Let $\u_{1\times n}$, $\X_{n\times n}$, and $\v_{n\times 1}$ be matrices of indeterminates,  $\Adj \X$ be the classical adjoint of $\X$, and $H(n)$ be the ideal $I_1(\u\X)+I_1(\X\v)+I_1(\v\u-\Adj \X)$. Vasconcelos has conjectured that $H(n)$ is a perfect Gorenstein ideal of grade $2n$. In this paper, we obtain the minimal free resolution of $H(n)$; and thereby establish Vasconcelos' conjecture. \endabstract
\title Ideals associated to two sequences and a matrix \endtitle 
\subjclass 13H10, 13D25 \endsubjclass
\leftheadtext{Andrew R. Kustin}
\rightheadtext{Two sequences and a matrix}
\author Andrew R. Kustin\endauthor
\address
Mathematics Department,
University of South Carolina,
Columbia, SC 29208\endaddress
\email kustin\@math.scarolina.edu \endemail
\endtopmatter

\document
 
\SectionNumber=\zero\tNumber=1
Let $\u_{1\times n}$, $\X_{n\times n}$, and $\v_{n\times 1}$ be matrices of indeterminates over a commutative noetherian ring $R_0$, and let   $H(n)$ be  the ideal 
$I_1(\u\X)+I_1(\X\v)+I_1(\v\u-\Adj \X)$ of the polynomial ring $R=R_0[\{u_i, v_i, x_{ij}\mid 1\le i,j\le n\}]$. 
 Vasconcelos has conjectured, in \cite{\rref{V91}, Conjecture 3.3.1}, that the ideal $H(n)$ is a perfect Gorenstein ideal of grade $2n$. In this paper, we obtain the minimal homogeneous resolution  of $R/H(n)$ by free $R-$modules; and thereby establish Vasconcelos' conjecture. 

In fact, we produce two resolutions of $R/H(n)$. The complex $\F$ of section \number\comF\ is never minimal, but it is relatively straightforward. The complex $\M$ of section \number\comM\ is a quotient of $\F$. It is more complicated than $\F$, but it is minimal. The exactness of $\F$ and $\M$ is established in section \number\xact. In section 5 we consider the singular locus and linkage history of $R/H(n)$; we also consider the  algebra structure of $\Tor_{\bullet}^R(R/H(n) , R_0 )$. At the end of the paper we record some open  questions  about the ring $R/H(n)$. 

 The arguments of sections \number\comF, \number\comM, and \number\xact\ are long, but routine. In section \number\comF\ we prove that $\F$ is a complex. In section \number\comM\ we split off a split exact subcomplex $\N$ of $\F$. In 
section \number\xact\ we apply the acyclicity lemma and reduce the problem  to one involving generic data with the parameter $n$ replaced by $n-1$. The most interesting part of the argument is the discovery of the complex  $\F$. This complex  is obtained by merging four Koszul complexes:$$ \matrix\format\c&\ \c&\ \c\\ \F(1) & \longleftrightarrow& \F(2)\\ \updownarrow&& \updownarrow\\ \F(3) & \longleftrightarrow& \F(4),\endmatrix\tag{*}$$where $\F(1)$ and $\F(4)$ are both Koszul complexes on the entries of $[\u\ \ \v]$, $\F(2)$ is the Koszul complex on the entries of $[\u\X\ \ \v]$, and $\F(3)$ is the Koszul complex on the entries of $[\u\ \ \X\v]$.
The arrows in (*) represent maps given by the various minors of $\X$. 

\bigpagebreak

\SectionNumber=\prelim\tNumber=1

\flushpar{\bf \number\SectionNumber.\quad Preliminary results.}

\medskip

In this paper ``ring'' means commutative noetherian ring with one.
The {\it grade} of a proper ideal $I$ in a ring $R$ is the length of the
longest regular sequence on $R$ in $I$. The ideal $I$ of $R$ is called {\it
perfect} if the grade of $I$ is equal to the projective dimension
 of the $R-$module $R/I$. The grade $g$ ideal $I$ is called {\it Gorenstein} if it is
perfect and $\Ext_{R}^{g}(R/I,R)\iso R/I.$ It follows from Bass \cite{\rref{Bss}, Proposition 5.1} that if $I$ is a Gorenstein ideal in a Gorenstein ring $R$, then $R/I$ is also a Gorenstein ring. 

Let $R$ be a ring. 
For any $R-$module $F$, we write $F^*=\Hom_R(F,R)$. If $f\:F\to G$ is a map of $R-$modules, we define $I_r(f)$ to be the image of the map $\W^rF\t(\W^rG)^*\to R$, which is induced by the map $\W^rf\:\W^rF\to\W^rG$. 

\definition{Definition \tnum{D1.2}} Let $R$ be a commutative ring. If $\u_{1\times n}$, $\X_{n\times n}$, and $\v_{n\times 1}$ are matrices with entries from $R$, then   $H(\u,\X,\v)$ is defined to be  the ideal 
$$I_1(\u\X)+I_1(\X\v)+I_1(\v\u-\Adj \X)$$ of $R$, where 
$\Adj \X$ is the classical adjoint of $\X$. (In other words,
 $\X\cdot \Adj \X$ and $\Adj \X\cdot \X$ are both equal to  $\det \X\cdot I$).
  \enddefinition

 Let $R$ be a commutative noetherian ring, and $F$ be 
 a free $R-$module of finite rank. We  make much use of the exterior algebras $\W^{\bullet}F$ and  $\W^{\bullet}F^*$. Each element of 
$F^{*}$  is a graded derivation on $\W^{\bullet}F$. In other words, $$\a_1\(\fakeht a_1^{[1]}\w\dots\w a_1^{[s]}\)=\sum_j (-1)^{j+1} \a_1(a_1^{[j]}) \cdot a_1^{[1]}\w\dots\w\widehat{a_1^{[j]}}\w\dots\w a_1^{[s]}\in \W^{j-1}F, $$ for all $\a_1\in F^*$ and $a_1^{[j]}\in F$. This action gives rise
to the 
 $\W^{\bullet}F^{*}-$module structure on $\W^{\bullet}F$. In particular, 
$$(\a_1\w \b_1)(a_s)= \a_1\(\fakeht \b_1(a_s)\),$$ for $\a_1, \b_1\in F^*$ and $a_s\in \W^sF$. 
The
$\W^{\bullet}F-$module 
structure  on $\W^{\bullet}F^{*}$ is obtained in an analogous manner. In
particular, if 
$a_{i}\in\W^{i}F$ and 
$\b_{j}\in \W^{j}F^{*}$, then 
$$a_{i}(\b_{j})\in \W^{j-i}F^{*}\qquad\text{and}\qquad 
\b_{j}(a_{i})\in \W^{i-j}F.$$ One   consequence of these two module structures is that $a_s(\a_s)=\a_s(a_s)\in R$ for all $a_s$ in $\W^sF$ and $\a_s\in \W^sF^*$. 
The following well known formulas show more of the  interaction between the
two module 
structures.
\proclaim{Proposition \tnum{A3}} Let $F$ be a free module over a commutative
noetherian ring $R$ and let $a,b\in \W^{\bullet}F$ and $\g\in\W^{\bullet}F^{*}$ be
homogeneous elements. 
\roster
\item"{(a)}" If $\deg a =1$, then 
$$ \(a(\g)\)(b)=a\w(\g(b))+(-1)^{1+\deg \g}\g(a\w b).$$ 
\item"{(b)}" If $\g \in \W^{\rank F}F^{*}$, then
$$\(a(\g)\)(b)=(-1)^{\nu}\(b(\g)\)(a), $$where $\nu=(\rank F -\deg a)(\rank F-
\deg b)$.\endroster\endproclaim
\remark{Note}The value for $\nu$ which is given above is correct. An incorrect value has appeared elsewhere in the literature.  \endremark

\proclaim{Corollary \tnum{31L3+}} Retain the hypotheses of Proposition \tref{A3}. If $b\in \W^{\rank F}F$, then $[a(\g)](b)$ is equal to $a\w\g(b)$. \endproclaim 

\demo{Proof} The proof is by induction on $\deg a$. The case $\deg a=1$ is established in Proposition \tref{A3}\,(a). If $a=a_1\w a'$, with  $\deg a_1=1$, then use the case $\deg a_1=1$ and the induction hypothesis to see that
$$\[\fakeht (a_1\w a')(\g)\](b)= \[a_1\(\fakeht a'(\g)\)\](b)= a_1 \w \(a'(\g)\fakeht\)(b) =(a_1\w a')\w \g(b). \qed$$\enddemo

\remark {Remark \tnum{R1}}Let  $F$ be a free module over a commutative ring $R$. The exterior algebra $\W^{\bullet}F$ comes equipped with  co-multiplication
$$\Delta\: \W^{\bullet}F\to \W^{\bullet}F\t\W^{\bullet}F.$$Co-multiplication is the algebra map which is induced by the diagonal map $F\to F\p F$. For example, 
if $a_1$, $a_1'$, and $a_1''$ are elements of $\W^1F$, then 
$$\Delta(a_1\w a_1'\w a_1'')= \left\{\matrix a_1\w a_1'\w a_1''\t 1 \in \W^3F\t \W^0F\\+\\ a_1\w a_1'\t a_1''- a_1\w a_1''\t a_1'+ a_1'\w a_1''\t a_1\in \W^2F\t \W^1F\\+\\ a_1\t a_1'\w a_1''-a_1'\t a_1\w a_1''+a_1''\t a_1\w a_1'\in \W^1F\t \W^2F\\+\\ 1\t a_1\w a_1'\w a_1'' \in \W^0F\t \W^3F\endmatrix\right.$$
Often, we will   use only one graded piece of the co-multiplication map. If $p+q=t$, then we write 
$$\Delta(a_t)=\sum\limits_j a_p^{[j]}\t a_{q}^{[j]}$$ to mean that the image of $a_t$ under the composition
$$\W^tF@>\text{inclusion}>> \W^{\bullet}F@>\Delta>>\W^{\bullet}F\t \W^{\bullet}F @>\text{projection}>> \W^pF\t \W^qF$$ is $\sum\limits_j a_p^{[j]}\w a_{q}^{[j]}.$ In particular, if $p=1$, $q=2$, and $a_3= a_1\w a_1'\w a_1''$, then 
$$\gather \Delta(a_3)=\sum_{j=1}^3 a_1^{[j]}\t a_2^{[j]},\quad\text{where}\\ a_1^{[1]}=a_1,\ a_2^{[1]}= a_1'\w a_1'',\ \ a_1^{[2]}=-a_1', \ a_2^{[2]}= a_1\w a_1'',\ \  a_1^{[3]}=a_1'',\ \text{and}\  a_2^{[3]}= a_1\w a_1'.\endgather $$
\endremark 

\proclaim{Lemma \tnum{31L2+}}Let  $F$ be a free module over a commutative ring $R$. Let $a_k$, $b_k$, and $c_k$  be elements of $\W^kF$ and  $\a_k$ and $\b_k$ be elements of $\W^kF^*$ for all integers $k$.  \roster
\item"{(a)}"  If $\Delta(\a_t)=\sum\limits_j\a_1^{[j]}\t\a_{t-1}^{[j]}$, then
$\sum\limits_j c_{t-2}(\a_{t-1}^{[j]})\w \a_1^{[j]}=2(-1)^{t-1}c_{t-2}(\a_t).$
\item"{(b)}"If $\Delta(\b_s)=\sum\limits_i \b_1^{[i]}\t \b_{s-1}^{[i]}$, then 
$\sum\limits_i\[\fakeht \b_1^{[i]}(a_2)\](\b_{s-1}^{[i]}) =-2a_2(\b_s).$
\endroster
\endproclaim

\demo{Proof}Apply $F$ to each side of (a). If $b_1$ is a fixed, but arbitrary, element of $F$, then 
$$\eightpoint \align & b_1\(\sum_j c_{t-2}(\a_{t-1}^{[j]})\w \a_1^{[j]}\) = \sum_j (b_1\w c_{t-2})(\a_{t-1}^{[j]})\cdot  \a_1^{[j]} -\sum_j b_1(\a_1^{[j]})\cdot c_{t-2}(\a_{t-1}^{[j]})\\&{ }
= (-1)^{t-1} (b_1\w c_{t-2})(\a_t) -c_{t-2}\(\fakeht b_1(\a_t)\)= (-1)^{t-1}2 b_1\(\fakeht c_{t-2}(\a_t) \). \endalign$$It suffices to prove (b) for $a_2=a_1\w a_1'$. In this case,
$$\align \sum_i\[\fakeht \b_1^{[i]}(a_1 \w a_1')\](\b_{s-1}^{[i]}) &{ }
= \sum_i  a_1 (\b_1^{[i]})\cdot a_1'(\b_{s-1}^{[i]}) -\sum_i  a_1' (\b_1^{[i]})\cdot a_1(\b_{s-1}^{[i]})\\
&{ }= a_1'\(\fakeht a_1(\b_s)\)-a_1\(\fakeht a_1'(\b_s)\)= -2 a_2(\b_s). \qed
\endalign$$
\enddemo

The following data is in effect throughout most of the paper. 

\definition{Data \tnum{SU}}Let $R$ be a commutative noetherian ring, $F$
be a free module of rank $n\ge 2$ over 
 $R$, $X\:F\to F^{*}$ be an $R-$module
homomorphism, and  $u$ and $v$ be elements of $F$. Fix  orientation
elements $\eta\in\W^{n}F$ and $\e_n\in \W^nF^*$, which are compatible in the sense that $e_n(\e_n)=1$.
\enddefinition

\remark{Note}We will always take $a_i$ and $b_i$ to be elements of $\W^iF$, and 
$\a_i$ and $\b_i$ to be elements of $\W^iF^*$.
\endremark

\remark{Remark \tnum{R1.8}}In the notation of Data \tref{SU}, let $\b_1^{[i]}\in \W^1F^*$ and $\a_t\in \W^tF^*$. In order to make sure that the reader understands our conventions, we give an expanded account of the symbol $\[\b_1^{[i]}\w (\W^{n-t}X^*)(\a_t[\eta])\](\eta)$, which appears in Definition \tref{36D1}: $\a_t[\eta]$ is the element of $\W^{n-t}F$ which is given by the module action of $\W^{\bullet}F^*$ on $\W^{\bullet}F$; the map  $(\W^{n-t}X^*)$ carries $\a_t[\eta]$ to an element of $\W^{n-t}F^*$; $\b_1^{[i]}\w (\W^{n-t}X^*)(\a_t[\eta])$ is an element of $\W^{n-t+1}F^*$; and the module action of $\W^{\bullet}F^*$ on $\W^{\bullet}F$ makes $\[\b_1^{[i]}\w (\W^{n-t}X^*)(\a_t[\eta])\](\eta)$ an element of $\W^{t-1}F$. \endremark 

 The following convention and calculation provide the connection between the coordinate free complexes $\F$ and $\M$ of sections \number\comF\ and \number\comM\ and the coordinate dependent ideals $H(\u,\X,\v)$ of Definition \tref{D1.2}.

\definition{Convention \tnum{conv}} Whenever we convert $u$, $X$, and $v$ from Data \tref{SU} into matrices $\u$, $\X$, and $\v$, we consider a pair of bases $e_1^{[1]},\dots,e_1^{[n]}$ for $F$ and $\e_1^{[1]},\dots,\e_1^{[n]}$ for $F^*$ which satisfy $\e_1^{[i]}(e_1^{[j]})=\delta_{ij}$ (the Kronecker delta), $e_1^{[1]}\w\dots\w e_1^{[n]}=\eta$, and $\e_1^{[n]}\w\dots\w \e_1^{[1]}=\e_n$. If $u=\sum\limits_i u_i e_1^{[i]}$, $X(e_1^{[j]})=\sum\limits_i x_{ij}\e_1^{[i]}$, and $v=\sum\limits_i v_i e_1^{[i]}$, then $\u=[u_1,\dots,u_n]$, $\X$ is the $n\times n$ matrix whose entry in row $i$ and column $j$ is $x_{ij}$, and $\v=\[\smallmatrix v_1\\ \vdots\\v_n\endsmallmatrix\]$. \enddefinition

\proclaim{Lemma \tnum{L3.10}} Adopt Data \tref{SU}. If $\u$, $\X$ and $\v$ are matrices which satisfy Convention \tref{conv}, then $$\[\e_1^{[j]}\w(\W^{n-1}X)(\e_1^{[i]}[\eta])\](\eta)=(-1)^{\frac{n(n-1)}{2}}(\Adj \X)_{ij}.$$\endproclaim
\demo{Proof} Let $\X(r_1,\dots,r_s;c_1,\dots,c_s)$ represent the determinant of the submatrix of $\X$ which consists of rows $r_1,\dots,r_s$ and columns $c_1,\dots,c_s$. The left side of the proposed identity is equal to 
$$\split & (-1)^{i+1}\[\e_1^{[j]}\w(\W^{n-1}X)\(e_1^{[1]}\w\dots\w\widehat{e_1^{[i]}}\w\dots\w
e_1^{[n]} \)\](\eta)\\ & { } = 
(-1)^{i+1}\[\e_1^{[j]}\w \(\sum\limits_k \X(1,\dots,\widehat{k},\dots,n;1,\dots,\widehat{i},\dots,n)  
\e_1^{[1]}\w\dots\w\widehat{\e_1^{[k]}}\w\dots\w \e_1^{[n]} \)\](\eta)\\ & { } =
(-1)^{i+j}\X(1,\dots,\widehat{j},\dots,n;1,\dots,\widehat{i},\dots,n)\cdot (\e_1^{[1]}\w\dots\w\e_1^{[n]})(\eta).\endsplit$$ The proof is complete since  $(\e_1^{[1]}\w\dots\w\e_1^{[n]})(\eta)=(-1)^{\frac{n(n-1)}{2}}$,   and  $$(-1)^{i+j}\X(1,\dots,\widehat{j},\dots,n;1,\dots,\widehat{i},\dots,n)=(\Adj\X)_{ij}. \qed$$\enddemo

Adopt Data \tref{SU}. If $a_1$ and  $b_{1}$ are elements of $F$, then the canonical identification of $F^{**}$ with $F$ yields 
$$\(\fakeht X(b_{1})\)(a_1)=b_1
\(\fakeht X^{*}(a_1)\).$$
Furthermore, the canonical identification of $(\W^jF)^*$ with $\W^jF^*$ gives 
$$\[(\W^{j}X)(b_{j})\](a_{j})
=b_j
\((\W^{j}X^{*})(a_{j})\)\tag\tnum{pg124}$$ for all $a_{j}$ and $b_{j}$ in $\W^{j}F$ and for all
$j$. Formula (\tref{pg124}) is the case ``$i=0$'' of part (a) of the following result.    

\proclaim{Observation \tnum{OBB}}Adopt Data \tref{SU}. Let $a_k$ and $b_k$ be elements of $\W^kF$ and  $\a_k$ and $\b_k$ be elements of $\W^kF^*$ for all integers $k$. \roster\item"{(a)}" $(\W^{i}X^{*})\[\((\W^{j}X)(b_{j})\)(a_{i+j})\]=
b_{j}\[\(\W^{i+j}X^{*}\)(a_{i+j})\].$
\item"{(b)}"   If $\Delta(a_t)=\sum\limits_j a_1^{[j]}\t a_{t-1}^{[j]}$ and $\Delta(\b_s)=\sum\limits_i\b_{1}^{[i]}\t\b_{s-1}^{[i]}$, then
$$\sum\limits_i\b_1^{[i]}(a_t)\t \b_{s-1}^{[i]}= \sum\limits_j  a_{t-1}^{[j]}\t a_1^{[j]}(\b_s). $$
\item"{(c)}" If $\Delta(\a_t)=\sum\limits_j \a_1^{[j]}\t \a_{t-1}^{[j]}$ and $\Delta(\b_s)=\sum\limits_i \b_1^{[i]}\t \b_{s-1}^{[i]}$,
 then 
$$\sum_i \( \[(\W^{n-1}X)[\b_1^{[i]}(\eta)]\](\eta)\)(\a_t)\t \b_{s-1}^{[i]}
=\sum_j\a_{t-1}^{[j]} \t \(\[(\W^{n-1}X^*)[\a_1^{[j]}(\eta)]\](\eta)\)(\b_s).
$$
\endroster 
\endproclaim 

\demo{Proof}  We apply each side of (a)   to the element
$a_{i}$ of $\W^{i}F$. Notice that $\((\W^{j}X)(b_{j})\)(a_{i+j})$ is in
$\W^{i}F$. Use 
(\tref{pg124}), together with the module actions of $\W^{\bullet}F$ and $\W^{\bullet}F^*$ on one another,   to see that 
$$\split & \((\W^{i}X^{*})  \[\((\W^{j}X)(b_{j})\)(a_{i+j})\]\fakeht
\)(a_{i})=\[\((\W^{j}X)(b_{j})\)(a_{i+j})\]\((\W^{i}X)(a_{i})\)\\ & \phantom{(\W^{i}X^{*})}  { } =
\[(\W^{i}X)(a_{i})\w 
(\W^{j}X)(b_{j})\](a_{i+j})  =\[(\W^{i+j}X)(a_{i}\w b_{j})\](a_{i+j})\\  &  \phantom{(\W^{i}X^{*})}{ }=
(a_{i}\w b_{j})\[(\W^{i+j}X^{*})(a_{i+j})\]
=\(b_{j}\[\(\W^{i+j}X^{*}\)(a_{i+j})\]\)(a_{i}).\endsplit$$   Both expressions in (b) are equal to  $   \sum\limits_i\sum\limits_j \b_1^{[i]}(a_1^{[j]})\cdot a_{t-1}^{[j]}\t \b_{s-1}^{[i]}$. Use the action of $F$ on $\W^{\bullet}F^*$, Proposition \tref{A3}\,(b), and (\tref{pg124}) to see that the left side of (c) is equal to 
$$\eightpoint  \align &  
\sum_i \sum_j \( \[(\W^{n-1}X)[\b_1^{[i]}(\eta)]\](\eta)\)(\a_1^{[j]})\cdot  \a_{t-1}^{[j]}\t\b_{s-1}^{[i]} \\& { }  =  (-1)^{n-1}
\sum_i \sum_j  \(\a_1^{[j]}[\eta]\) \[(\W^{n-1}X)[\b_1^{[i]}(\eta)]\]\cdot  \a_{t-1}^{[j]}\t\b_{s-1}^{[i]}\\&{ } =(-1)^{n-1} \sum_i \sum_j [\b_1^{[i]}(\eta)]\[(\W^{n-1}X^*)(\a_1^{[j]}(\eta))\]\cdot  \a_{t-1}^{[j]}\t\b_{s-1}^{[i]} \\& { } = \sum_j\a_{t-1}^{[j]} \t \sum_i \(\[(\W^{n-1}X^*)\(\a_1^{[j]}(\eta)\)\](\eta)\)(\b_1^{[i]})\cdot\b_{s-1}^{[i]},
\endalign $$
which is equal to the right side of (c). 
  $\qed$ 
\enddemo

\bigpagebreak

\SectionNumber=\comF\tNumber=1

\flushpar{\bf \number\SectionNumber.\quad The complex $\F$.}

\medskip

 The modules and maps which comprise the complex $\F$ are given in Definition \tref{36D1}. Recall the conventions  of Remarks \tref{R1} and \tref{R1.8}. If $t$ is an integer with $t\le -1$ or $n+1\le t$, then $\W^tF=0$. 

\definition {Definition \tnum{36D1}}Adopt Data \tref{SU}. The module $\F_r$ of $\F$ is defined to be  $$\F_r= \F_r(1)\p\F_r(2)\p\F_r(3)\p \F_r(4),\quad\text{where} \quad
\F_r(i)=\sum_{t\in\Z}\F_r(i)^{(t)}\quad\text{for}$$  
$$\matrix \format\l&\quad\l \\ \F_r(1)^{(t)}= \W^tF^*\t \W^{r+1-t}F^*,&
\F_r(2)^{(t)}=\W^tF\t \W^{r-t}F^*,\\ \vspace{5pt}
\F_r(3)^{(t)}=\W^tF\t \W^{r-t}F^*,\ \text{and} &
\F_r(4)^{(t)}= \W^tF^*\t \W^{r-1-t}F^*.\endmatrix$$
The maps $$f_r(i)^{(t)}\: \F_r(i)^{(t)} @>>> \F_{r-1}$$ are given below. 

\flushpar If $s+t=r+1$, $\Delta(\a_t)=\sum\limits_j\a_1^{[j]}\t\a_{t-1}^{[j]}$, and
$\Delta(\b_s)=\sum\limits_i\b_{1}^{[i]}\t\b_{s-1}^{[i]}$, then
$$ f_r(1)^{(t)}(\a_t\t\b_s)= \left\{\matrix 
\a_t\t v(\b_s)\in \F_{r-1}(1)^{(t)} \\+\\
(-1)^r u(\a_t)\t \b_s \in \F_{r-1}(1)^{(t-1)}\\ +\\\sum\limits_i\[\b_1^{[i]}\w (\W^{n-t}X^*)(\a_t[\eta])\](\eta)\t\b_{s-1}^{[i]}\in \F_{r-1}(2)^{(t-1)}\\+\\
 \sum\limits_j \[\a_1^{[j]}\w(\W^{n-s}X)(\b_s[\eta])\](\eta)\t\a_{t-1}^{[j]}\in \F_{r-1}(3)^{(s-1)}
 .\endmatrix \right.$$
If $s+t=r$ and $\Delta(a_t)=\sum\limits_j a_1^{[j]}\t a_{t-1}^{[j]}$, then
$$f_r(2)^{(t)}(a_t\t \b_s) = \left\{\matrix
(-1)^r 
(\W^tX)(a_t)\t \b_s \in \F_{r-1}(1)^{(t)}\\+\\
a_t\t v(\b_s)\in \F_{r-1}(2)^{(t)}\\+\\
(-1)^{r+1}[X^*(u)](a_t)\t\b_s\in \F_{r-1}(2)^{(t-1)}\\+\\  
 \sum\limits_j (\W^{t-1}X)(a_{t-1}^{[j]})\t a_1^{[j]}(\b_s)\in \F_{r-1}(4)^{(t-1)}.
\endmatrix \right.$$
If $s+t=r$ and
$\Delta(a_t)=\sum\limits_j a_{1}^{[j]}\t a_{t-1}^{[j]}$, then $$f_r(3)^{(t)}(a_t\t \b_s) = \left\{\matrix
(-1)^{r+1} \b_s\t (\W^tX^*)(a_t)\in \F_{r-1}(1)^{(s)}\\+\\
(-1)^{r+1} a_t\t u(\b_s)  \in \F_{r-1}(3)^{(t)}\\+\\ [X(v)](a_t) \t  \b_s\in \F_{r-1}(3)^{(t-1)}\\+\\
- \sum\limits_j a_1^{[j]}(\b_s)\t (\W^{t-1}X^*)(a_{t-1}^{[j]})\in \F_{r-1}(4)^{(s-1)}.
\endmatrix \right.$$
If $s+t=r-1$, then $$f_r(4)^{(t)}(\a_t\t \b_s) = \left\{\matrix
(-1)^r \[(\W^{n-t}X^*)(\a_t[\eta])\](\eta)\t \b_s \in \F_{r-1}(2)^{(t)}\\+\\ 
(-1)^r  \[(\W^{n-s}X)(\b_s[\eta])\](\eta)\t\a_t\in \F_{r-1}(3)^{(s)}\\+\\
\a_t\t v(\b_s)\in \F_{r-1}(4)^{(t)}\\+\\ 
(-1)^r u(\a_t)\t \b_s\in \F_{r-1}(4)^{(t-1)} .
\endmatrix \right.$$
\enddefinition

\proclaim{Proposition \tnum{P1}}
The modules and maps of Definition \tref{36D1} form a complex 
$$ (\F,f):\quad  0@>>> \F_{2n+1}@>>>\dots @>>>\F_r @>f_r >> \F_{r-1} @>>> \dots @>>> \F_{-1}@>>>0.$$
\endproclaim
\remark{Note} When we want to emphasize the data which was used to construct $(\F,f)$, we write $\F\,[u,X,v]$.
\endremark 
\demo{Proof}We prove that $f_{r-1}\circ f_r(\ell )^{(t)}=0$  for $1\le \ell \le 4$. In each case   we write  $$\eightpoint f_{r-1}\circ f_r(\ell)^{(t)}(y_t\t \b_s)= A+B+C+D,\ A=\sum_{k=1}^4 A(k),\ B=\sum_{k=1}^4 B(k),\ C=\sum_{k=1}^4 C(k),\ \text{and}\ D=\sum_{k=1}^4 D(k),$$ where $y_t=\a_t$ if $\ell= 1$ or $4$, and   $y_t=a_t$ if $\ell= 2$ or $3$.  

\flushpar{\bf The case $\ell=1$.} Let  $s+t=r+1$, $\Delta(\a_t)=\sum\limits_j\a_1^{[j]}\t\a_{t-1}^{[j]}$, and
$\Delta(\b_s)=\sum\limits_i\b_{1}^{[i]}\t\b_{s-1}^{[i]}$.  We have $$\eightpoint \spreadlines{5pt} \allowdisplaybreaks \align &  A=f_{r-1}(1)^{(t)} 
\(\a_t\t v(\b_s)\fakeht\), \\
&B= (-1)^r  f_{r-1}(1)^{(t-1)}\(\fakeht u(\a_t)\t \b_s\), \\
&C=f_{r-1}(2)^{(t-1)}\(\sum\limits_i\[\b_1^{[i]}\w (\W^{n-t}X^*)(\a_t[\eta])\](\eta)\t\b_{s-1}^{[i]}\), \\
& D= f_{r-1}(3)^{(s-1)}\(\sum\limits_j \[\a_1^{[j]}\w(\W^{n-s}X)(\b_s[\eta])\](\eta)\t\a_{t-1}^{[j]}\), \endalign  $$
 $$  \eightpoint \spreadlines{5pt} \allowdisplaybreaks \align &A(1)= \a_t\t v(v(\b_s)) \in \F_{r-2}(1)^{(t)}, \\&
A(2)= (-1)^{r-1} u(\a_t) \t v(\b_s) \in \F_{r-2}(1)^{(t-1)}, \\&
A(3)= - \sum\limits_i\[\b_1^{[i]}\w (\W^{n-t}X^*)(\a_t[\eta])\](\eta)\t v(\b_{s-1}^{[i]}) \in \F_{r-2}(2)^{(t-1)},\\&
A(4)= \sum\limits_j \[\a_1^{[j]}\w(\W^{n-s+1}X)\([v(\b_s)][\eta]\fakeht\)\](\eta)\t\a_{t-1}^{[j]}\in \F_{r-2}(3)^{(s-2)},\\&
B(1)= (-1)^r u(\a_t)\t v(\b_s) \in \F_{r-2}(1)^{(t-1)},\\&
B(2)= - u(u(\a_t)) \t \b_s \in \F_{r-2}(1)^{(t-2)},\\&
B(3)= (-1)^r \sum\limits_i\[\b_1^{[i]}\w (\W^{n-t+1}X^*)\([u(\a_t)][\eta]\fakeht\)\](\eta)\t\b_{s-1}^{[i]}\in \F_{r-2}(2)^{(t-2)},\\&
B(4)= (-1)^{r+1} \sum\limits_j \[\a_1^{[j]}\w(\W^{n-s}X)(\b_s[\eta])\](\eta)\t u(\a_{t-1}^{[j]})\in \F_{r-2}(3)^{(s-1)},\\&
C(1)= (-1)^{r-1} \sum\limits_i (\W^{t-1}X) \(\[\b_1^{[i]}\w (\W^{n-t}X^*)(\a_t[\eta])\](\eta)\) \t\b_{s-1}^{[i]}\in \F_{r-2}(1)^{(t-1)},\\&
C(2)= \sum\limits_i\[\b_1^{[i]}\w (\W^{n-t}X^*)(\a_t[\eta])\](\eta)\t v(\b_{s-1}^{[i]}) \in \F_{r-2}(2)^{(t-1)},\\&
C(3)= (-1)^r \sum\limits_i [X^*(u)] \(\[\b_1^{[i]}\w (\W^{n-t}X^*)(\a_t[\eta])\](\eta)\)\t\b_{s-1}^{[i]}\in \F_{r-2}(2)^{(t-2)},\\&
C(4)=\text{\ the 
$\F_{r-2}(4)^{(t-2)}-$component of $f_{r-1}(2)^{(t-1)}\(\sum\limits_i\[\b_1^{[i]}\w (\W^{n-t}X^*)(\a_t[\eta])\](\eta)\t\b_{s-1}^{[i]}\)$} ,\\&
D(1)= (-1)^r \sum\limits_j \a_{t-1}^{[j]} \t (\W^{s-1}X^*) \( \[\a_1^{[j]}\w(\W^{n-s}X)(\b_s[\eta])\](\eta)\)  \in \F_{r-2}(1)^{(t-1)},\\&
D(2)= (-1)^r \sum\limits_j \[\a_1^{[j]}\w(\W^{n-s}X)(\b_s[\eta])\](\eta)\t u(\a_{t-1}^{[j]})\in \F_{r-2}(3)^{(s-1)},\\&
D(3)= \sum\limits_j [X(v)] \(\[\a_1^{[j]}\w(\W^{n-s}X)(\b_s[\eta])\](\eta)\) \t\a_{t-1}^{[j]}\in \F_{r-2}(3)^{(s-2)},\quad \text{and}\\&
D(4)= \text{\  the 
$\F_{r-2}(4)^{(t-2)}-$component of $ f_{r-1}(3)^{(s-1)}\(\sum\limits_j \[\a_1^{[j]}\w(\W^{n-s}X)(\b_s[\eta])\](\eta)\t\a_{t-1}^{[j]}\)$}.\endalign$$ 
Observe that $$0=A(1)=B(2)=A(2)+B(1)=A(3)+C(2)=A(4)+D(3)=B(3)+C(3)=B(4)+D(2).$$
Use the module action of $\W^{\bullet}F^*$ on $\W^{\bullet}F$, together with Observation \tref{OBB}\,(a) and Proposition \tref{A3}\,(b), to see that   $$\allowdisplaybreaks \align C(1)= { } &   
(-1)^{n-t} (-1)^{r-1} 
\sum\limits_i (\a_t[\eta])\( (\W^{n-1}X)\(\b_1^{[i]}[\eta]\) \)\t\b_{s-1}^{[i]}
\\ { }= { } &  (-1)^{r-1}   \sum\limits_i \(\[(\W^{n-1}X)\(\b_1^{[i]}[\eta]\)\](\eta)\)(\a_t) \t \b_{s-1}^{[i]}, \ \text{and}\\  D(1)=  { } &  (-1)^{n-s}(-1)^r 
\sum\limits_j\a_{t-1}^{[j]}\t (\b_s[\eta]) \((\W^{n-1}X^*)(\a_1^{[j]}[\eta] )\)\\ { }= { } & 
(-1)^r 
\sum\limits_j\a_{t-1}^{[j]}\t \(\[(\W^{n-1}X^*)(\a_1^{[j]}[\eta] )\](\eta)\)(\b_s)
.\endalign $$ Apply Observation \tref{OBB}\,(c) to see that $C(1)+D(1)=0$.  

We prove $C(4)+D(4)=0$ by showing that $(c_{t-2}\t 1)*(C(4)+D(4))=0$ for all $c_{t-2}\in \W^{t-2}F$, where
$$(c_{t-2}\t 1)*(\a_{t-2}\t \b)=c_{t-2}(\a_{t-2})\cdot \b\in \W^{\bullet}F^*,\tag\tnum{T2.6}$$
for $\a_{t-2}\t \b\in \W^{t-2}F^*\t \W^{\bullet}F^*$.
Let $c_{t-2}$ be a fixed, but arbitrary, element of $\W^{t-2}F$. 
Observe that $$\eightpoint \spreadlines{5pt} \align & (c_{t-2}\t 1)*\(\ \text{ the 
$\F_{r-2}(4)^{(t-2)}-$component of $f_{r-1}(2)^{(t-1)}(a_{t-1}\t \b_{s-1})$}\)\\&{ } =  \sum\limits_j c_{t-2}\((\W^{t-2}X)(a_{t-2}^{[j]})\)\cdot a_1^{[j]}(\b_{s-1}) = 
  \sum\limits_j \((\W^{t-2}X^*)(c_{t-2})\)\(
a_{t-2}^{[j]}\)\cdot a_1^{[j]}(\b_{s-1})\\&{ } = (-1)^t  \[\Fakeht \[(\W^{t-2}X^*)(c_{t-2})\](a_{t-1})\](\b_{s-1}),\endalign $$where $\Delta(a_{t-1})=\sum\limits_j a_1^{[j]}\t a_{t-2}^{[j]}$.  It follows that   $$\eightpoint \align (c_{t-2}\t 1)*C(4)   = { } &
(-1)^t  
\sum_i \[\Fakeht[(\W^{t-2}X^*)(c_{t-2})]\(\[\b_1^{[i]}\w(\W^{n-t}X^*)(\a_t(\eta))\](\eta)\)\](\b_{s-1}^{[i]})\\ { } = { } &
  \sum_i \(\Fakeht\[\b_1^{[i]}\w(\W^{n-2}X^*)\(\fakeht 
c_{t-2}\w \a_t(\eta)\)\](\eta)\)(\b_{s-1}^{[i]}).\endalign $$ Use Lemma \tref{31L2+}\,(b) and Corollary \tref{31L3+}
  to see that  
$$\eightpoint \align (c_{t-2}\t 1)*C(4) &{ } =-2   
\(\[(\W^{n-2}X^*)\(\fakeht 
c_{t-2}\w \a_t(\eta)\)\](\eta)\Fakeht\)(\b_s)\\&{ } =-2   
\(\[(\W^{n-2}X^*)\(\fakeht 
[c_{t-2}(\a_t)](\eta)\)\](\eta)\Fakeht\)(\b_s).\endalign $$
 Observe, also, that  
$$\eightpoint\spreadlines{5pt} \align &  (c_{t-2}\t 1) *\(\fakeht\ \text{the $\F_{r-2}(4)^{(t-2)}-$component of 
$f_{r-1}(3)^{(s-1)}(a_{s-1}\t \b_{t-1})$}\)\\ &{ } =- \sum\limits_j c_{t-2}\(a_1^{[j]}(\b_{t-1})\)\cdot (\W^{s-2}X^*)(a_{s-2}^{[j]})  = (-1)^{t-1}    (\W^{s-2}X^*) \(\sum\limits_j [c_{t-2}(\b_{t-1})](a_1^{[j]})\cdot a_{s-2}^{[j]}\)\\ &{ }= (-1)^{t-1}   (\W^{s-2}X^*)\(\fakeht [c_{t-2}(\b_{t-1})](a_{s-1})\),\endalign$$ where $\Delta(a_{s-1})=\sum\limits_j a_1^{[j]}\t a_{s-2}^{[j]}$.
It follows  that 
$(c_{t-2}\t 1)*D(4) $ is equal to $$ (-1)^{t-1} (\W^{s-2}X^*)
\(\sum\limits_j \[c_{t-2}(\a_{t-1}^{[j]})\]\(\[\a_1^{[j]}\w (\W^{n-s}X) (\b_s(\eta))\](\eta) \)\).$$ Apply Lemma \tref{31L2+}\,(a), Observation \tref{OBB}\,(a), and Proposition \tref{A3}\,(b) to see that 
$$\nopagebreak \eightpoint\align & (c_{t-2}\t 1)*D(4)=2   (\W^{s-2}X^*)
\(\Fakeht \[c_{t-2}(\a_{t})\w (\W^{n-s}X) (\b_s[\eta])\](\eta) \)\\&{ }=
2 (\b_s[\eta]) \[(\W^{n-2}X^*)\(\fakeht [c_{t-2}(\a_t)](\eta)\)\] \\&{} =2 \[\Fakeht \[(\W^{n-2}X^*)\(\fakeht[c_{t-2}(\a_t)](\eta)\)\](\eta)\](\b_s)=-(c_{t-2}\t 1)*C(4);\endalign $$ therefore, $C(4)+D(4)=0$ and $f_{r-1}\circ f_r(1)^{(t)}(\a_t\t\b_s)=0$. 
 
\flushpar{\bf The case $\ell=2$.} Let  $s+t=r$, $\Delta(a_t)=\sum\limits_j a_1^{[j]}\t a_{t-1}^{[j]}$, and $\Delta(\b_s)=\sum\limits_i\b_{1}^{[i]}\t\b_{s-1}^{[i]}$.  We have  $$\eightpoint \spreadlines{5pt} \allowdisplaybreaks \align & 
A=(-1)^r 
f_{r-1}(1)^{(t)}\( (\W^tX)(a_t)\t \b_s \)
 ,\\
&B= f_{r-1}(2)^{(t)}\( a_t\t v(\b_s)\fakeht\),\\
&C=(-1)^{r+1}f_{r-1}(2)^{(t-1)}\(\fakeht [X^*(u)](a_t)\t\b_s\),\\
&D= f_{r-1}(4)^{(t-1)}\( \sum\limits_j (\W^{t-1}X)(a_{t-1}^{[j]})\t a_1^{[j]}(\b_s)\), 
\endalign $$  
$$ \eightpoint \spreadlines{5pt} \allowdisplaybreaks \align 
&A(1)= (-1)^r (\W^tX)(a_t)\t v(\b_s)
\in \F_{r-2}(1)^{(t)},\\&
A(2)=- u\((\W^tX)(a_t)\)\t \b_s \in \F_{r-2}(1)^{(t-1)},\\&
A(3)=(-1)^r  \sum\limits_i\[\b_1^{[i]}\w (\W^{n-t}X^*)\(\fakeht [(\W^tX)(a_t)][\eta]\)\](\eta)\t\b_{s-1}^{[i]} \in \F_{r-2}(2)^{(t-1)},\\&
A(4)= (-1)^r \sum\limits_j \[X(a_1^{[j]})\w(\W^{n-s}X)(\b_s[\eta])\](\eta)\t (\W^{t-1}X)(a_{t-1}^{[j]})   \in \F_{r-2}(3)^{(s-1)},\\&
B(1)=    (-1)^{r-1} (\W^tX)(a_t)\t v(\b_s)  \in \F_{r-2}(1)^{(t)},\\&
B(2)=  a_t\t v[v(\b_s) ]  \in \F_{r-2}(2)^{(t)},\\&
B(3)=   (-1)^{r}[X^*(u)](a_t)\t v(\b_s) \in \F_{r-2}(2)^{(t-1)},\\&
B(4)=  \sum\limits_j (\W^{t-1}X)(a_{t-1}^{[j]})\t a_1^{[j]}(v[\b_s])  \in \F_{r-2}(4)^{(t-1)},\\&
C(1)=   (\W^{t-1}X)\(\fakeht [X^*(u)](a_t)\)\t \b_s 
\in \F_{r-2}(1)^{(t-1)},\\&
C(2)=(-1)^{r+1}[X^*(u)](a_t)\t v(\b_s)
 \in \F_{r-2}(2)^{(t-1)},\\&
C(3)=-[X^*(u)]\([X^*(u)](a_t)\)\t\b_s  \in \F_{r-2}(2)^{(t-2)},\\&
C(4)=(-1)^{r}\sum\limits_j u\((\W^{t-1}X)(a_{t-1}^{[j]})\)\t a_1^{[j]}(\b_s) \in \F_{r-2}(4)^{(t-2)},\\&
D(1)= (-1)^{r-1} \sum\limits_j  \[(\W^{n-t+1}X^*) \[ \( (\W^{t-1}X)(a_{t-1}^{[j]}) \) (e_n)\] \Fakeht \](e_n)  \t a_1^{[j]}(\b_s)  \in \F_{r-2}(2)^{(t-1)},\\&
D(2)= (-1)^{r-1}  \sum\limits_j    
\[(\W^{n-s+1}X)\[ \(a_1^{[j]}(\b_s)\) (e_n)\] \Fakeht\] (e_n)
\t (\W^{t-1}X)(a_{t-1}^{[j]}) \in \F_{r-2}(3)^{(s-1)},\\&
D(3)=\sum\limits_j (\W^{t-1}X)(a_{t-1}^{[j]})\t v\(a_1^{[j]}(\b_s)\) \in \F_{r-2}(4)^{(t-1)},\quad\text{and}\\&
D(4)= (-1)^{r-1}\sum\limits_j u\((\W^{t-1}X)(a_{t-1}^{[j]})\)\t a_1^{[j]}(\b_s)\in \F_{r-2}(4)^{(t-2)}.\endalign$$ 
Observe that $$\eightpoint 0=B(2)=C(3)=A(1)+B(1)=A(2)+C(1)=A(4)+D(2)=B(3)+C(2)=B(4)+D(3)=C(4)+D(4).$$
Furthermore, Observation \tref{OBB}\,(a) and (b), Proposition \tref{A3}\,(b), and Corollary \tref{31L3+}  give 
$$ \eightpoint \matrix\format\l\\ A(3) =   (-1)^r [(\W^nX)(\eta)](\eta)\cdot \sum\limits_i\b_1^{[i]}(a_t)\t \b_{s-1}^{[i]}\\   = (-1)^r[(\W^nX)(\eta)](\eta)\cdot \sum\limits_j  a_{t-1}^{[j]}\t a_1^{[j]}(\b_s)=-D(1); \endmatrix  \tag\tnum{T2.9}$$thus, $f_{r-1}\circ f_r(2)^{(t)}(a_t\t \b_s)=0$. 
 
\flushpar{\bf The case $\ell=3$.}  Let  $s+t=r$, $\Delta(a_t)=\sum\limits_j a_1^{[j]}\t a_{t-1}^{[j]}$, and $\Delta(\b_s)=\sum\limits_i\b_{1}^{[i]}\t\b_{s-1}^{[i]}$.  We have $$\eightpoint \spreadlines{5pt} \allowdisplaybreaks \align &
A= (-1)^{r+1} f_{r-1}(1)^{(s)}\( \b_s\t (\W^tX^*)(a_t)\)
,\\
&B= (-1)^{r+1} f_{r-1}(3)^{(t)}\(\fakeht  a_t\t u(\b_s) \)
,\\
&C= f_{r-1}(3)^{(t-1)}\(\fakeht [X(v)](a_t) \t  \b_s\)
, \\
&D= - f_{r-1}(4)^{(s-1)}\(\sum\limits_j a_1^{[j]}(\b_s)\t (\W^{t-1}X^*)(a_{t-1}^{[j]})\)
, \endalign $$
$$ \eightpoint  \allowdisplaybreaks \spreadlines{5pt} \align 
&A(1)=(-1)^{r+1}\b_s\t v\((\W^tX^*)(a_t)\)
 \in \F_{r-2}(1)^{(s)},\\&
A(2)=  u(\b_s) \t (\W^tX^*)(a_t)  \in \F_{r-2}(1)^{(s-1)},\\&
A(3)=(-1)^{r+1}\sum\limits_j \[X^*(a_1^{[j]})\w(\W^{n-s}X^*)(\b_s[\eta])\](\eta)\t
(\W^{t-1}X^*)( a_{t-1}^{[j]} ) \in \F_{r-2}(2)^{(s-1)},\\&
A(4)=(-1)^{r+1}\sum\limits_i\[\b_1^{[i]}\w (\W^{n-t}X) \[ \fakeht [(\W^tX^*)(a_t)] [\eta]\]
 \Fakeht \](\eta)\t 
  \b_{s-1}^{[i]} \in \F_{r-2}(3)^{(t-1)},\\&
B(1)= - u(\b_s)\t (\W^tX^*)(a_t)\in \F_{r-2}(1)^{(s-1)},\\&
B(2)=- a_t\t u(u(\b_s))   \in \F_{r-2}(3)^{(t)},\\&
B(3)= (-1)^{r+1}[X(v)](a_t) \t  u(\b_s )\in \F_{r-2}(3)^{(t-1)},\\&
B(4)=(-1)^{r}\sum\limits_j a_1^{[j]}(u[\b_s])\t (\W^{t-1}X^*)(a_{t-1}^{[j]}) \in \F_{r-2}(4)^{(s-2)},\\&
C(1)=(-1)^{r} \b_s\t (\W^{t-1}X^*)\([X(v)](a_t) \fakeht\) \in \F_{r-2}(1)^{(s)},\\&
C(2)= (-1)^{r} [X(v)](a_t)\t u(\b_s)    \in \F_{r-2}(3)^{(t-1)},\\&
C(3)= [X(v)]\([X(v)](a_t)\fakeht \) \t  \b_s \in \F_{r-2}(3)^{(t-2)},\\&
C(4)= \sum\limits_j a_1^{[j]}(\b_s)\t v \[(\W^{t-1}X^*)(a_{t-1}^{[j]})\]\in \F_{r-2}(4)^{(s-1)},\\&
D(1)= (-1)^{r} \sum\limits_j   \[(\W^{n-s+1}X^*) \[\(a_1^{[j]}(\b_s)\)(\eta)\]\Fakeht\] (\eta)
  \t (\W^{t-1}X^*)(a_{t-1}^{[j]})\in \F_{r-2}(2)^{(s-1)},\\&
D(2)=(-1)^{r} \sum\limits_j \[(\W^{n-t+1}X) \Fakeht \[ \((\W^{t-1}X^*)(a_{t-1}^{[j]})
\)(e_n) \]   \](e_n) \t a_1^{[j]}(\b_s)\in \F_{r-2}(3)^{(t-1)},\\&
D(3)=- \sum\limits_j a_1^{[j]}(\b_s)\t v \[(\W^{t-1}X^*)(a_{t-1}^{[j]})\]\in \F_{r-2}(4)^{(s-1)},\quad\text{and}\\&
D(4)=(-1)^{r} \sum\limits_j u\[ a_1^{[j]}(\b_s)\] \t (\W^{t-1}X^*)(a_{t-1}^{[j]}) \in \F_{r-2}(4)^{(s-2)}.\endalign $$
Observe that
$$\eightpoint 0=B(2)=C(3)=A(1)+C(1)= A(2)+B(1)= A(3)+D(1)=B(3)+C(2)=B(4)+D(4)
=C(4)+D(3).$$
The argument of (\tref{T2.9}) gives  
$ A(4)+D(2)=0$; and therefore, $f_{r-1}\circ f_r(3)^{(t)}(a_t\t \b_s)=0$.
 
\flushpar{\bf The case $\ell=4$.}  Let $s+t=r-1$. We have  $$\eightpoint \spreadlines{5pt} \allowdisplaybreaks \align & 
A= (-1)^r f_{r-1}(2)^{(t)}\( \[(\W^{n-t}X^*)(\a_t[\eta])\](\eta)\t \b_s \)
,\\
&B= (-1)^r f_{r-1}(3)^{(s)}\( \[(\W^{n-s}X)(\b_s[\eta])\](\eta)\t\a_t\),\\
&C=f_{r-1}(4)^{(t)}\(\fakeht \a_t\t v(\b_s)\)
 , \\
&D= (-1)^r f_{r-1}(4)^{(t-1)}\( \fakeht u(\a_t)\t \b_s\),\endalign $$   
$$\eightpoint  \spreadmatrixlines{5pt}\allowdisplaybreaks  \align  
&A(1)=- (\W^tX) \( \[(\W^{n-t}X^*)(\a_t[\eta])\](\eta)\) \t \b_s\in \F_{r-2}(1)^{(t)},\\&
A(2)=(-1)^r\[(\W^{n-t}X^*)(\a_t[\eta])\](\eta)\t v(\b_s) \in \F_{r-2}(2)^{(t)},\\&
A(3)= [X^*(u)] \(\[(\W^{n-t}X^*)(\a_t[\eta])\](\eta)\) \t \b_s \in \F_{r-2}(2)^{(t-1)},\\&
A(4)=\text{\  the $\F_{r-2}(4)^{(t-1)}-$component of $(-1)^r    f_{r-1}(2)^{(t)}\(\fakeht \[(\W^{n-t}X^*)(\a_t[\eta])\](\eta)\t \b_s
\)$},\\&
B(1)= \a_t \t  (\W^sX^*)\( \[(\W^{n-s}X)(\b_s[\eta])\](\eta)\) \in \F_{r-2}(1)^{(t)},\\&
B(2)=  \[(\W^{n-s}X)(\b_s[\eta])\](\eta)\t u(\a_t) \in \F_{r-2}(3)^{(s)},\\&
B(3)=(-1)^r [X(v)]\( \[(\W^{n-s}X)(\b_s[\eta])\](\eta)\) \t\a_t\in \F_{r-2}(3)^{(s-1)},\\&
B(4)=\text{\  the $\F_{r-2}(4)^{(t-1)}-$component of $(-1)^r    f_{r-1}(3)^{(s)}
 \(\fakeht 
\[(\W^{n-s}X)(\b_s[\eta])\](\eta)\t \a_t 
\)$},\\&
C(1)= (-1)^{r-1} \[(\W^{n-t}X^*)(\a_t[\eta])\](\eta)\t v(\b_s )\in \F_{r-2}(2)^{(t)},\\&
C(2)= (-1)^{r-1}  \[(\W^{n-s+1}X)((v[\b_s])[\eta])\](\eta)\t\a_t\in \F_{r-2}(3)^{(s-1)},\\&
C(3)= \a_t\t v(v(\b_s) )\in \F_{r-2}(4)^{(t)},\\&
C(4)=(-1)^{r-1}  u(\a_t)\t v(\b_s)  \in \F_{r-2}(4)^{(t-1)},\\&
D(1)= - \[(\W^{n-t+1}X^*)\([u(\a_t)] [\eta]\fakeht\)\](\eta)\t \b_s \in \F_{r-2}(2)^{(t-1)},\\&
D(2)=-  \[(\W^{n-s}X)(\b_s[\eta])\](\eta)\t u(\a_t) \in \F_{r-2}(3)^{(s)},\\&
D(3)=(-1)^r u(\a_t)\t v(\b_s) \in \F_{r-2}(4)^{(t-1)},\quad\text{and}\\&
D(4)=- u(u(\a_t))\t \b_s\in \F_{r-2}(4)^{(t-2)}.\endalign $$
Observe that $$\eightpoint 0=C(3)=D(4)=A(1)+B(1)=A(2)+C(1)= A(3)+D(1)=B(2)+D(2)= B(3)+C(2)=C(4)+D(3).$$
Let $c_{t-1}$ be a fixed, but arbitrary,  element of $\W^{t-1}F$. Employ the trick of (\tref{T2.6}).  We see that  $$ \eightpoint \split   & (c_{t-1}\t 1)* (\text { the  $\F_{r-2}(4)^{(t-1)}-$component
of $f_{r-1}(2)^{(t)}(a_t\t \b_s) $} )\\ & { } =
  \sum\limits_j c_{t-1}\( (\W^{t-1}X)(a_{t-1}^{[j]})\) \cdot a_1^{[j]}(\b_s){ } = 
   \(\sum\limits_j  \[(\W^{t-1}X^*)(c_{t-1})\] (a_{t-1}^{[j]}) \cdot a_1^{[j]}\)(\b_s)\\ & { }= 
(-1)^{t-1}    \( \[(\W^{t-1}X^*)(c_{t-1})\] (a_{t}) \)(\b_s),\endsplit $$where $\Delta(a_{t})=\sum\limits_j a_1^{[j]}\t a_{t-1}^{[j]}$; and therefore, it follows that  
$(c_{t-1}\t 1) *A(4)$ is equal to $$\eightpoint \split & (-1)^{t-1+r}    \( \[(\W^{t-1}X^*)(c_{t-1})\] \(
\[(\W^{n-t}X^*)(\a_t[\eta])\](\eta)\) \Fakeht \)(\b_s) \\
&{ }= (-1)^{t-1+r}     \( \[(\W^{n-1}X^*) \([c_{t-1}(\a_t)][\eta]\) \](\eta) \Fakeht \)(\b_s)\\
 &{ } =(-1)^{n-s } (-1)^{t-1+r}    [\b_s(\eta)] \[\Fakeht (\W^{n-1}X^*) \(\fakeht [c_{t-1}(\a_t)] (\eta)\)\].\endsplit$$
We also see that  $(c_{t-1}\t 1)* \( \text{ the 
$\F_{r-2}(4)^{(t-1)}-$component of 
$f_{r-1}(3)^{(s)}( b_{s}\t \a_t)$}\)$ is equal to 
$$\split  - \sum\limits_i c_{t-1}
 \( b_1^{[i]}(\a_t)\) \cdot (\W^{s-1}X^*)(b_{s-1}^{[i]})
 = { } &(-1)^{t}  (\W^{s-1}X^*)\[\Fakeht
\sum\limits_i \(\fakeht c_{t-1}
 (\a_t)\) (b_1^{[i]})\cdot  b_{s-1}^{[i]}\]\\
{ }  = { } &(-1)^{t}  (\W^{s-1}X^*) \[\(\fakeht c_{t-1}
 (\a_t)\) (b_s)\Fakeht\],\endsplit $$where $\Delta(b_{s})=\sum\limits_i b_1^{[i]}\t b_{s-1}^{[i]}$;
therefore,  $  (c_{t-1}\t 1)*B(4)$ is equal to $$\eightpoint \split &
(-1)^{t+r}  (\W^{s-1}X^*) \[\(\fakeht c_{t-1}
 (\a_t)\) \( \[(\W^{n-s}X)(\b_s[\eta])\](\eta) \Fakeht\)\] \\ { } ={ } & (-1)^{t+r}(-1)^{n-s}   (\b_s[\eta])\[ (\W^{n-1}X^*)\(\fakeht [c_{t-1}(\a_t)][\eta]\)\]
= -(c_{t-1}\t 1)*A(4).\endsplit $$
It follows that 
$A(4)+B(4)=0$ and the proof is complete. 
$\qed$\enddemo

\remark{Remark \tnum{R2.4}} Suppose that the data  of \tref{SU} is graded. Let  $X$ be  a homogeneous homomorphism of degree 1, and let $u$ and $v$ be homogeneous elements of $F$ of degree $d_u$ and $d_v$, respectively. If $d_u +d_v=n-1$, then it is easy to check that $\F$ is a graded complex with  homogeneous maps of degree zero, provided the   grading on $\F$ is   given by: 
$$\split & \F_r(1)^{(t)}=R^{\binom{n}{t}\binom{n}{r+1-t}}\[-\(\fakeht tn-t+d_v(r+1-2t)\)\] \\ &
 \F_r(2)^{(t)}=R^{\binom{n}{t}\binom{n}{r-t}}\[-\(\fakeht tn+d_v(r-2t)\)\] \\ &
 \F_r(3)^{(t)}=R^{\binom{n}{t}\binom{n}{r-t}}\[-\(\fakeht (r-t)n+(d_v+1)(2t-r)\)\] \\ &
 \F_r(4)^{(t)}=R^{\binom{n}{t}\binom{n}{r-1-t}}\[-\(\fakeht (t+1)n-t +d_v(r-1-2t)\)\].\endsplit$$ 
\endremark 

\bigpagebreak

\SectionNumber=\comM\tNumber=1

\flushpar{\bf \number\SectionNumber.\quad The complex $\M$\,.}

\medskip

Theorem \tref{T3.9} is the main result in this section. Its proof appears after the proof of  Proposition \tref{P3.16}.   When we want to emphasize the data which was used to construct $(\M,m)$, we write $\M\,[u,X,v]$.  The case $n=2$ is handled  in Proposition \tref{R3.21}.
 
\proclaim{Theorem \tnum{T3.9}} Adopt Data \tref{SU} with $3\le n$. Let  $(\M,m)$ be the maps and modules of  Definitions \tref{42D1} and \tref{42D1'}\,$($c$)$. The following statements hold.
\roster
\item"{(a)}"The maps and modules of $(\M,m)$ form a complex
$$\M:\quad 0@>>> \M\,_{2n}@>>>\dots @>>> \M\,_r@>m_r>>\M\,_{r-1}@>>>\dots @>>> \M\,_0.$$
\item"{(b)}"  Let   $\u$, $\X$ and $\v$ be matrices which satisfy Convention \tref{conv}. If 
  $H$ is  the ideal $H(\u,\X,(-1)^{\frac{n(n-1)}{2}}\v)$ of Definition \tref{D1.2}, then the homology  $H_0(\M)$ is equal to $R/H$,
\item"{(c)}"If $\F$ is the complex   of Definition \tref{36D1}, then $H_r(\F)=H_r(\M)$ for all $r$.  
\item"{(d)}" Each map $m_r$ of  $\M$ satisfies $I_1(m_r)\subseteq I_1(u)+I_1(v)+I_1(X)$.\endroster
\endproclaim 

In section \number\xact\ we prove that $\M$ is acyclic   whenever Data \tref{SU} is sufficiently generic (in the sense of Corollary \tref{C5.6}). If, in addition, the data is    local or graded (in the sense of Remark \tref{R2.4}), then assertion (d) of the above result ensures that $\M$ is a minimal resolution.
 Some notation must be fixed before we can   describe the modules of $\M$. 

\definition{Definition \tnum{D3.4}}Adopt Data \tref{SU}. For each integer $s$, let $$\mu_s\:F\t \W^sF^*\to \W^{s-1}F^*\quad \text{and}\quad \bx_s\: \W^sF^* \to \W^{n-1}F\t \W^{s-1}F^* $$ be the homomorphisms which are given by 
$$ \mu_s(a_1\t\a_s)=a_1(\a_s)\quad\text{and}\quad \bx_s(\a_s ) =\sum_i 
\a_1^{[i]}(e_n)\t \a_{s-1}^{[i]},$$ where $\Delta(\a_s)=\sum\limits_i \a_1^{[i]} \t \a_{s-1}^{[i]}$.
 \enddefinition

  \remark{Observation} Notice that $\mu_{s+1}$ is a surjection for all $s$, except $s=n$; and $\bx_{s+1}$ is a split injection for all $s$, except $s=-1$.
\endremark 

 \definition{Definition \tnum{D3.5}}Retain the notation of Definition \tref{D3.4}. For each integer $s$, define homomorphisms $$\ell_s\:\W^sF^*\to F\t \W^{s+1}F^*\quad\text{and}\quad \lambda_{s}\:\W^{n-1}F\t \W^{s}F^* \to \W^{s+1}F^*,\quad \text{by}$$
\roster
\item"{(a)}"$\ell_s$ is a fixed splitting of $\mu_{s+1}$ for  $s\neq n$;  
 \item"{(b)}"$\ell_n=0$;
\item"{(c)}" $\lambda_s$ is a fixed splitting of $\bx_{s+1}$ for $s\neq -1$;
 and 
\item"{(d)}"$\lambda_{-1}=0$. 
\endroster
\enddefinition

\remark{Remark \tnum{41D0}}The maps $\ell_s$ and $\lambda_s$ have been chosen so that 
$$\split \mu_{s+1}\circ \ell_s=\id &\quad\text{for all integers $s$, except $s=n$, and }\\ 
 \lambda_{s}\circ \bx_{s+1}=\id &\quad\text{for all $s$,   except $s=-1$.}\endsplit$$ 
%Furthermore, it is obvious   that 
%$$\ell_s\circ \mu_{s+1}\circ \ell_s=\ell_s \quad\text{for all integers $s$.} $$
\endremark

 \definition{Definition \tnum{D3.6}}Retain the notation of Definition \tref{D3.5} with $3\le n$. 
For $i=2$ and $3$, define submodules $\[\F_r(i)^{(n-1)}\]'$ and $\[\F_r(i)^{(n-1)}\]''$ of $\F_r(i)^{(n-1)}$, and submodules 
$\[\F_r(i)^{(1)}\]'$ and $\[\F_r(i)^{(1)}\]''$ of $\F_r(i)^{(1)}$ by 
$$\allowdisplaybreaks \split \[\F_r(i)^{(n-1)}\]'= { } &\Ker\[\F_r(i)^{(n-1)}= \W^{n-1}F\t \W^{r+1-n}F^* @>\lambda_{r+1-n}>> \W^{r+2-n}F^*\], \\
\[\F_r(i)^{(n-1)}\]''= { } &\Im \[\W^{r+2-n}F^* @>\bx_{r+2-n} >> \W^{n-1}F\t \W^{r+1-n}F^* =\F_r(i)^{(n-1)}\],\\
\[\F_r(i)^{(1)}\]'= { } &\Ker \[ \F_r(i)^{(1)}= F\t\W^{r-1}F^*
@>\mu_{r-1}>> \W^{r-2}F^*\],\ \text{and}\\
\[\F_r(i)^{(1)}\]''= { } &\Im \[ \W^{r-2}F^*@>\ell_{r-2}>> F\t\W^{r-1}F^*= \F_r(i)^{(1)}\].\endsplit $$
\enddefinition
The following statements are immediate consequences of Definition \tref{D3.6}. 
\proclaim{Observation \tnum{O3.8}} If  $i=2$ or $3$,  then 
\roster
\item"{(a)}"$\F_r(i)^{(1)}=\[\F_r(i)^{(1)}\]'\p\[\F_r(i)^{(1)}\]''$ for all $r$,
\item"{(b)}"$\F_r(i)^{(n-1)}=\[\F_r(i)^{(n-1)}\]'\p\[\F_r(i)^{(n-1)}\]''$ for all $r$,
\item"{(c)}" $\[\F_r(i)^{(1)}\]''=0$,  for $ r\le 1$,
\item"{(d)}"$\[\F_r(i)^{(n-1)}\]''=0$, for $2n-1\le r$. 
\item"{(e)}"$\[\F_r(i)^{(1)}\]'=0$,  for $n+1\le r$, and
\item"{(f)}" $\[\F_r(i)^{(n-1)}\]'=0$,  for $ r\le n-1$.
\endroster\endproclaim
    
\definition{Definition \tnum{42D1}}Adopt Data \tref{SU} with $3\le n$. The module $\M\,_r$ of $\M=\M\,[u,X,v]$ is obtained as follows.  Let $\widehat{\M\,}_r$ represent the following
submodule of $\F_r$:
$$\eightpoint \split \widehat{\M}\,_r=  { } &
\sum_{t\notin\{0,n,r+1,r+1-n\}}\F_r(1)^{(t)}+ \[\F_r(2)^{(n-1)}\]' +\sum_{2\le t\le n-2}\F_r(2)^{(t)}+ \[\F_r(2)^{(1)}\]' \\\vspace{5pt}{ } & { }
+\[\F_r(3)^{(n-1)}\]'+\sum_{2\le t\le n-2 }\F_r(3)^{(t)} +\[\F_r(3)^{(1)}\]'+\sum_{t\notin\{0,n,r-1,r-1-n\}}\F_r(4)^{(t)}. \endsplit$$The submodule $\M\,_r$ of $\F_r$ is defined by  $$ \M\,_r= \cases \widehat{\M}\,_0+\F_0(3)^{(0)}, & \text{if $r=0$,}\\\widehat{\M}\,_2+\[\F_2(3)^{(1)}\]''
,& \text{if $r=2$,}\\\widehat{\M}\,_{2n-2}+\[\F_{2n-2}(3)^{(n-1)}\]'',& \text{if $r=2n-2$,}\\\widehat{\M}\,_{2n}+\F_{2n}(3)^{(n)}
,& \text{if $r=2n$, and}\\
\widehat{\M}\,_{r}
,& \text{for all other $r$.}
\endcases$$ \enddefinition

\remark{Remark \tnum{R3.1}} Adopt   the grading hypotheses of   Remark \tref{R2.4}. If $n=3$ and $d_u=d_v=1$, then $\M$ is   
$$\eightpoint 0\to R(-9)\to R(-7)^{15}\to  R(-6)^{35}\to \matrix R(-4)^{21}\\\p\\R(-5)^{21}\endmatrix \to R(-3)^{35}\to R(-2)^{15}\to R.$$
If $n=4$, $d_u=1$, and $d_v=2$, then $\M$ is   
$$\eightpoint \split 0\to R(-16)\to \matrix R(-13)^{20}\\\p\\ R(-14)^4\endmatrix \to \matrix R(-10)^6\\\p\\ R(-11)^{24}\\\p\\ R(-12)^{61}\endmatrix \to \matrix R(-9)^{56}\\\p\\ R(-10)^{80}\\\p\\R(-11)^{36}\endmatrix \to \matrix R(-6)^{10}\\\p\\ R(-7)^{24}\\\p\\ R(-8)^{140}\\\p\\ R(-9)^{24}\\\p\\ R(-10)^{10}\endmatrix \to \matrix R(-5)^{36}\\\p\\ R(-6)^{80}\\\p\\ R(-7)^{56}\endmatrix \to \matrix R(-4)^{61}\\\p\\ R(-5)^{24}\\\p\\ R(-6)^{6}\endmatrix  &\\ \to \matrix R(-2)^4\\\p\\ R(-3)^{20}\endmatrix \to R.&\endsplit $$
 In general, if $4\le n$, then   $\M\,_0=R$, 
$$\eightpoint  \allowdisplaybreaks\spreadlines{5pt}\align   \M\,_1= { } & { }R^{n^2}\[-(n-1)\fakeht\]\p R^n\[-(1+d_u)\fakeht\]\p
 R^n\[-(1+d_v)\fakeht\],\\\allowdisplaybreak
\M\,_2={ } & { } R^{\binom{n}{2}n}\[\fakeht -(n-1+d_u)\] \p R^{\binom{n}{2}n}\[\fakeht -(n-1+d_v)\] \p R^{\binom{n}{2}}\[\fakeht -(2+2d_u)\] \p R^{2n^2-1}\[\fakeht -n\] \\ & { }\p
R^{\binom{n}{2}}\[\fakeht -(2+2d_v)\],\\ 
\M\,_r ={ } & { } \sum\limits_{t=1}^{r} R^{\binom{n}{t}\binom{n}{r+1-t}}\[-\(\fakeht tn-t+d_v(r+1-2t)\)\]
 \p
R^{n\binom{n}{r-1}-\binom{n}{r-2} }\[-\(\fakeht n+d_v(r-2)\)\] 
\\ & { }  \p \sum\limits_{t=2}^{r}
R^{\binom{n}{t}\binom{n}{r-t}}\[-\(\fakeht tn+d_v(r-2t)\)\]
 \p 
R^{n\binom{n}{r-1}-\binom{n}{r-2} }\[-\(\fakeht 
(r-1)n+(d_v+1)(2-r)\)\]
\\ & { } \p \sum\limits_{t=2}^{r} 
R^{\binom{n}{t}\binom{n}{r-t}}\[-\(\fakeht (r-t)n+(d_v+1)(2t-r)\)\] 
\\ & { }  \p  
\sum\limits_{t=1}^{r-2}  R^{\binom{n}{t}\binom{n}{r-1-t}}\[-\(\fakeht (t+1)n-t +d_v(r-1-2t)\)\]
\ \text{for $3\le r\le n-2$},\\
\M\,_{n-1} ={ } & { } 
\sum\limits_{t=1}^{{n-1}} 
R^{\binom{n}{t}^2}\[-\(\fakeht tn-t+d_v(n-2t)\)\]
 \p
R^{n\binom{n}{2}-\binom{n}{3}}\[-\(\fakeht n+d_v(n-3)\)\]
\\ & { }\p   \sum\limits_{t=2}^{{n-2}}  
R^{\binom{n}{t}\binom{n}{t+1}}\[-\(\fakeht tn+d_v(n-1-2t)\)\]
\p R^{n\binom{n}{2}-\binom{n}{3}}\[-\(\fakeht n+d_u(n-3)\)\]
 \\ & { }\p \sum\limits_{t=2}^{{n-2}}
 R^{\binom{n}{t}\binom{n}{t+1}}\[-\(\fakeht (n-1-t)n+(d_v+1)(2t-n+1)\)\]  
\\ & { } \p   
\sum\limits_{t=1}^{n-3}
 R^{\binom{n}{t}\binom{n}{t+2}}\[-\(\fakeht (t+1)n-t +d_v(n-2-2t)\)\],\ \text{and}\\
\M\,_{n}={ } & { } \sum\limits_{t=2}^{n-1} R^{\binom{n}{t}\binom{n}{t-1}}\[-\(\fakeht tn-t+d_v(n+1-2t)\)\] \p  R^{n^2-\binom{n}{2}}\[-\(\fakeht n+d_v(n-2)\)\]
\\ & { }\p \sum\limits_{t=2}^{{n-2}} 
 R^{\binom{n}{t}^2}\[-\(\fakeht tn+d_v(n-2t)\)\]
\p
R^{n^2-\binom{n}{2}}\[-\(\fakeht 2n-2+d_u(n-2)\)\]
 \\ & { }
 \p   R^{n^2-\binom{n}{2}}\[-\(\fakeht n+d_u(n-2)\)\]
\p \sum\limits_{t=2}^{{n-2}}  R^{\binom{n}{t}^2}\[-\(\fakeht (n-t)n+(d_v+1)(2t-n)\)\]
\\ & { }\p R^{n^2-\binom{n}{2}} \[-\(\fakeht  2n-2+d_v(n-2)\)\]
\\ & { }\p \sum\limits_{t=1}^{{n}-2} R^{\binom{n}{t}\binom{n}{t+1}}\[-\(\fakeht (t+1)n-t +d_v(n-1-2t)\)\].
\endalign $$Furthermore, if $\M\,_r=\sum\limits_i R^{b_i}[-m_i]$, then $\M\,_{2n-r}= \sum\limits_iR^{b_i}[-(n^2-m_i)]$. 
\endremark

 \definition{Convention \tnum{T3.13}}  For each  statement  ``S'', let $$\chi(\text{S})= \cases 1,&\text{if S is true, and}\\ 0,&\text{if S is false.}\endcases$$In particular, $\chi(i=j)$ has the same value as the Kronecker delta $\delta_{ij}$.\enddefinition 

\definition{Definition \tnum{D3.12'}}Adopt Data \tref{SU} with $3\le n$. For each integer $r$,  let $\widehat{\N}\,_r$  be the following
submodule of $\F_r$:
 $$ \split \widehat{\N}\,_r=  { } & \chi(n\le r\le 2n-1) \cdot \F_r(1)^{(n)}+ 
\chi(n\le r\le 2n-2) \cdot \F_r(1)^{(r+1-n)}
\\ { } & { } +\chi(0\le r\le n) \cdot
\F_r(2)^{(0)}+ \chi(2\le r\le n+1)\[\F_r(2)^{(1)}\]''\\ { } & { } +\chi(1\le r\le n) \cdot \F_r(3)^{(0)}+\chi(3\le r\le n+1)  \cdot \[\F_r(3)^{(1)}\]''\\ { } & { }+\chi(n+1\le r\le 2n+1)  \cdot
\F_r(4)^{(n)}+\chi(n+1\le r\le 2n)  \cdot\F_r(4)^{(r-1-n)}.\endsplit$$ Let $(\N,n)$ be the subcomplex of $(\F,f)$ which is given by $$\N\,_r=\widehat{\N}\,_r + f_{r+1}(\widehat{\N}\,_{r+1})\quad\text{and}\quad n_r=f_r|_{\N\,_r}.$$ For each integer $r$, let 
$\L\,_r$   be the following
submodule of $\F_r$: $$  \split \L\,_r=   { } &  \chi(-1\le r\le n-1)  \cdot \F_r(1)^{(0)}
+\chi(0\le r \le n-1) \cdot   \F_r(1)^{(r+1)}  \\ { } & { }+\chi(n-1\le r\le 2n-2) \cdot \[\F_r(2)^{(n-1)}\]'' +\chi(n\le r\le 2n) \cdot \F_r(2)^{(n)}
\\ { } & { }
+\chi(n-1\le r\le 2n-3)  \cdot \[\F_r(3)^{(n-1)}\]''+\chi(n\le r\le 2n-1) \cdot \F_r(3)^{(n)}\\ { } & { }+\chi(1\le r\le n) \cdot \F_r(4)^{(0)} +\chi(2\le r \le n)  \cdot \F_r(4)^{(r-1)}.  \endsplit
$$
\enddefinition
 
\remark{Remark \tnum{R3.12''}}Use Observation \tref{O3.8} in order to see that 
  $$\F_r=\L\,_r\p\M\,_r\p\widehat{\N}\,_r\quad\text{ for all $r$}.$$   This  decomposition gives rise to projection maps 
$$\pi^{\L}_r\:\F_r\to\L\,_r,\quad \pi^{\M}_r\:\F_r\to\M\,_r,\quad\text{and}\quad
\pi^{\Nh}_r\:\F_r\to\widehat{\N}\,_r.$$For example, $\pi^{\L}_r$ is the map which annihilates $\M\,_r\p\widehat{\N}\,_r$, but restricts to give the identity on $\L\,_r$. \endremark 
 
\definition{Definition \tnum{42D1'}}Retain the notation of Definition \tref{D3.12'}. 
\roster\item"{(a)}" For each integer $r$, define 
$\tau_r\:\L\,_r\to \Nh\,_{r+1}$ by
 $$\eightpoint \matrix\format \l&\ \l\\  \tau_r(1)^{(0)}(1\t \b_{r+1})=(-1)^{r+1} 1\t \b_{r+1} \in \F_{r+1}(2)^{(0)}, &\text{for $-1\le r\le n-1 $,}\\
\tau_r(1)^{(r+1)}(\a_{r+1}\t 1)=(-1)^r 1\t \a_{r+1} \in \F_{r+1}(3)^{(0)}, &\text{for $0\le r\le n-1$,}\\
\tau_r(2)^{(n-1)} (\sigma_{r+2-n}(\b_{r+2-n}))=  \left\{\matrix 
\e_n\t\b_{r+2-n} \in \F_{r+1}(1)^{(n)}\\+\\
\d_{r\, n-1}(-1)^{n}v( \b_{1})\t\e_n \in \F_n(3)^{(0)},\endmatrix \right.&\text{for   $n-1\le r\le 2n-2$,}\\
\tau_r(2)^{(n)} (e_n\t\b_{r-n})= \left\{\matrix 
(-1)^{r+1}   \e_n \t \b_{r-n}\in \F_{r+1}(4)^{(n)}\\+\\
\d_{r\, n} \b_0\cdot  \ell_{n-1}[u(\e_n)]\in \[\F_{n+1}(3)^{(1)}\]'',\endmatrix \right.
&\text{for  $n\le r\le 2n$,}\\
\tau_r(3)^{(n-1)} (\sigma_{r+2-n}(\b_{r+2-n}))=\left\{\matrix 
\b_{r+2-n} \t \e_n \in \F_{r+1}(1)^{(r+2-n)}\\+\\
-\d_{r\,n-1}  u(\b_{1}) \t \e_n \in \F_n(2)^{(0)},
\endmatrix \right.&\text{for $ n-1\le r\le 2n-3$,}\\
\tau_r(3)^{(n)} (e_n\t\b_{r-n})=\left \{\matrix (-1)^{r+1} \b_{r-n}\t   \e_n \in \F_{r+1}(4)^{(r-n)}\\+\\ \d_{r\,n}(-1)^{n} \b_0\cdot  \ell_{n-1}[v(\e_n)]\in\[\F_{n+1}(2)^{(1)}\]''\endmatrix\right.
 &\text{for $n\le r\le 2n-1$,}\\
\tau_r(4)^{(0)}(1\t \b_{r-1})= \ell_{r-1}(\b_{r-1})\in \[\F_{r+1}(2)^{(1)}\]'',&\text{for $1\le r\le n$, and} \\
\tau_r(4)^{(r-1)}(\a_{r-1}\t1)=-\ell_{r-1}(\a_{r-1})\in \[\F_{r+1}(3)^{(1)}\]'',&\text{for $2\le r\le n$.} \endmatrix $$
\item"{(b)}"For each integer $r$, define $\psi_r\:\F_r\to\M\,_r$ by 
$$\psi_r|_{\M\,_r}=\id,\quad \psi_r|_{\Nh\,_r}=0\quad\text{and}\quad \psi_r|_{\L\,_r}=-\pi^{\M}_r\circ f_{r+1}\circ \tau_r.$$
\item"{(c)}" For each integer $r$, define $m_r\:\M\,_r\to\M\,_{r-1}$ to be the composition $$\M\,_r@>\incl>>\F_r@>f_r>>\F_{r-1}@>\psi_{r-1}>>\M\,_{r-1}.$$
\item"{(d)}" For each integer $r$, define $\rho_r\:\M\,_r\to\F_{r}$ by
$$\rho_r=\incl_r-\tau_{r-1}\circ\pi_{r-1}^{\L}\circ f_r.$$
\endroster
\enddefinition

\remark{Note} The definition of $\tau_r(i)^{(n-1)}$, for $i=2$ or $3$, is legitimate because Remark \tref{41D0} guarantees that $\lambda_{r+1-n}\circ\sigma_{r+2-n}=\id$, provided $n-1\le r$.
%%In the definition of $\tau_r(i)^{(n-1)}$, for $i=2$ or $3$, recall, from %Remark \tref{41D0}, that $\lambda_{r+1-n}\circ\sigma_{r+2-n}=\id$, provided %$n-1\le r$. 
\endremark

The technical part of the proof of Theorem \tref{T3.9} is contained in the proof of the next result. 

\proclaim{Lemma \tnum{L3.15}} In the notation of Definition \tref{42D1'}, the maps
$$\L\,_r@>\tau_r>> \Nh\,_{r+1}\quad\text{and}\quad \Nh\,_{r+1} @>\incl>> \F_{r+1}@>f_{r+1}>>\F_r @>\pi^{\L}_r >> \L\,_r$$ are inverses of one another. \endproclaim

\demo{Proof}We show that $\pi_r^{\L}\circ
f_{r+1}\circ \tau_r (x)=x$ for all $x\in \L\,_r$ and $\tau_r\circ \pi^{\L}_r\circ 
f_{{r+1}}(y)=y$ for all $y\in \Nh\,_{r+1}.$ There are eight cases. We first
 fix an integer $r$, with $-1\le r\le n-1$.  Let  $x=1\t
\b_{{r+1}}\in\F_{r}(1)^{(0)}$ and $y=1\t
\b_{{r+1}}\in\F_{r+1}(2)^{(0)}$. Observe that $$\tau_r(x)=(-1)^{r+1}y\quad \text{and}\quad \pi^{\L}_r\circ 
f_{{r+1}}(y) = (-1)^{r+1} x. $$
In the second case, we take $0\le r \le n-1$. Let $x=\a_{r+1}\t 1\in \F_r(1)^{(r+1)}$ and
 $y$ equal $ 1\t\a_{r+1}\in \F_{r+1}(3)^{(0)}$. Observe that 
$$ \tau_r(x)=(-1)^r y\quad  \text{and}\quad   \pi^{\L}_r\circ f_{r+1}(y)= (-1)^rx.$$In the third case, we have   $n-1\le r\le 2n-2$. Let $x=
\sigma_{r+2-n}(\b_{r+2-n})\in\[\F_r(2)^{(n-1)}\]'' $ and
 $y=\e_n\t \b_{r+2-n}\in \F_{r+1}(1)^{(n)}$. If $x'=\d_{r\, n-1}\cdot \e_n\t v(\b_{1})\in\F_{n-1}(1)^{(n)}$ and $y'$ is equal to $ \d_{r\, n-1}\cdot v(\b_{1})\t \e_n\in  \F_n(3)^{(0)}$, then   
$$ \eightpoint \tau_r (x)= y+(-1)^ny',\quad \tau_r(x')= (-1)^{n-1}y',\quad \pi^{\L}_r\circ f_{r+1}(y)=x+x',\quad \text{and}\quad \pi^{\L}_r\circ f_{r+1}(y')=(-1)^{n-1} x'.$$
In the fourth case, we consider $n\le r\le 2n$. Let $x=e_n\t\b_{r-n}\in  \F_r(2)^{(n)}$ and $y$ equal $
\e_n\t\b_{r-n}\in\F_{r+1}(4)^{(n)}$.  If $x'=\d_{r\,n}\cdot u(\e_n)\t\b_{0}\in \F_n(4)^{(n-1)}$ and $y'=\d_{r\, n}\cdot \b_{0}\cdot \ell_{n-1}[u(\e_n)]$ in $\[\F_{n+1}(3)^{(1)}\]''$, then   
$$\eightpoint   \tau_r (x)=(-1)^{r+1}y+y',\quad  \tau_r (x')=- y',\quad 
\pi^{\L}_r\circ f_{r+1}(y)=(-1)^{r+1}x+(-1)^{r+1}x',\quad \text{and}\ 
\pi^{\L}_r\circ f_{r+1}(y')=-x'.
$$
In case five, we have $n-1\le r\le 2n-3$. Let $x=\sigma_{r+2-n}(\b_{r+2-n})\in
\[\F_r(3)^{(n-1)}\]''$ and $y=\b_{r+2-n}\t\e_n\in\F_{r+1}(1)^{(r+2-n)}$.
If $x'=\d_{r\, n-1}\cdot u(\b_1)\t \e_n\in \F_{n-1}(1)^{(0)}$ and $y'=\d_{r\, n-1}\cdot u(\b_1)\t \e_n\in \F_n(2)^{(0)}$, then
$$\eightpoint \tau_r (x)=y-y',\quad  \tau_r (x')=(-1)^n y',\quad 
\pi^{\L}_r\circ f_{r+1}(y)=x+(-1)^{n}x',\quad \text{and}\quad 
\pi^{\L}_r\circ f_{r+1}(y')=(-1)^n x'.
$$ 
In the sixth case, we consider 
 $n\le r\le 2n-1$. Let $x=e_n\t\b_{r-n}\in  \F_r(3)^{(n)}$ and $y=\b_{r-n}\t\e_n\in \F_{r+1}(4)^{(r-n)}$. 
If $x'=\d_{r\,n}\cdot \b_0\t v(\e_n)\in \F_n(4)^{(0)}$ and $y'$ is equal to $\d_{r\, n}\cdot \b_0\cdot \ell_{n-1}(v[\e_n])\in \[\F_{n+1}(2)^{(1)}\]''$, then 
$$\eightpoint \tau_r (x)=(-1)^{r+1}y+(-1)^ny',\quad  \tau_r (x')=  y',\quad 
\pi^{\L}_r\circ f_{r+1}(y)=(-1)^{r+1}x+ x',\quad \text{and}\quad 
\pi^{\L}_r\circ f_{r+1}(y')= x'.
$$ 
In case seven, we have
  $1\le r\le n$. Let $x=1\t\b_{r-1}\in\F_r(4)^{(0)}$ and $y=\ell_{r-1}(\b_{r-1})$ in $\[\F_{r+1}(2)^{(1)}\]''$. Observe that 
$ \tau_r (x)= y$  and $ \pi^{\L}_r\circ f_{r+1}(y)= x$.  Finally, we take  $2\le r \le n$. Let $x=\a_{r-1}\t 1\in  \F_r(4)^{(r-1)}$ and
$y=  
\ell_{r-1}(\a_{r-1})\in\[\F_{r+1}(3)^{(1)}\]''$.
The proof is complete because $ \tau_r (x)= -y$  and $ \pi^{\L}_r\circ f_{r+1}(y)= -x$.
 $\qed$\enddemo

\proclaim{Proposition \tnum{P3.16}} Adopt Data \tref{SU} with $3\le n$. Let 
$(\F,f)$ be the complex  of  Definition \tref{36D1} and $(\N,n)$ be the subcomplex of $\F$ of Definition \tref{D3.12'}.  
\roster
\item"{(a)}" The complex $(\N,n)$   is split exact.
\item"{(b)}" The modules  and maps   $\{m_r\:\M\,_r\to\M\,_{r-1}\}$ of Definitions  \tref{42D1} and \tref{42D1'}\,$($c$)$ form a complex, which we denote $(\M,m)$.
\item"{(c)}" The maps $\{\psi_r\:\F_r\to \M\,_r\}$ of Definition \tref{42D1'}\,$($b$)$ form a map of complexes; furthermore,   
$$0\to (\N,n)@>\incl>> (\F,f) @>\psi>> (\M,m)\to 0$$is a short exact sequence of complexes. 
\item"{(d)}" The maps $\{\rho_r\:\M\,_r\to \F_r\}$ of Definition \tref{42D1'}\,$($d$)$ form a map of complexes; furthermore,  the composition
$$\M\,_r@>\rho_r>> \F_r@>\psi_r>>\M\,_r$$ is the identity map.
\endroster
\endproclaim

\demo{Proof}  Let  
$s_r\:\F_r\to \F_{r+1}$ be the map which is given by 
 $$ s_r|_{\L\,_r}=\tau_r,\quad s_r|_{\M\,_r}=0, \quad \text{and}\quad
s_r|_{\Nh\,_r}=0.$$ It follows, from Lemma \tref{L3.15}, that $$s_r|_{\Nh\,_r}=0\quad\text{and}\quad 
s_{r-1}\circ f_r|_{\Nh\,_r}=\id|_{\Nh\,_r}.$$ 
Assertion (a) is established    because the maps $s_r|_{\N\,_r}$   form a homotopy on $\N$ in the sense that 
$$s_{r-1}|_{\N\,_{r-1}}\circ n_r+n_{r+1}\circ s_r|_{\N\,_r}=\id|_{\N\,_r}\ \text{for all $r$.}$$  
We next show that $$\M\,_r+\N\,_r=\F_r\quad\text{for all $r$.}\tag\tnum{T3.18}$$ If $x_r\in\L\,_r$, then Lemma \tref{L3.15} gives $x_r=\pi_r^{\L}\circ f_{r+1}\circ \tau_r(x_r)$; and therefore, $x_r-f_{r+1}(\tau_r(x_r))$ is in $\M\,_r+\Nh\,_r$. It follows that $$\L\,_r\subseteq \M\,_r+\Nh\,_r+f_{r+1}(\Nh\,_{r+1})=\M\,_r+\N\,_r,$$ and (\tref{T3.18}) is established by Remark \tref{R3.12''}. 
Now we prove (b). 
Observe that $$m_r=\pi_{r-1}^{\M}(1-f_r\circ\tau_{r-1}\circ \pi_{r-1}^{\L})\circ f_r\quad\text{and}\quad m_{r+1}=\pi_r^{\M}\circ f_{r+1}(1-\tau_r\circ \pi_r^{\L}\circ f_{r+1}).$$ We know, from Remark \tref{R3.12''},  that $$\id|_{\F_r}= \pi_r^{\L}+\pi_r^{\M}+\pi_r^{\Nh}\quad\text{ and}\quad f_r\circ f_{r+1}=0;\tag\tnum{T3.40}$$ therefore, we see that
$$f_r\circ m_{r+1}= f_r\circ \pi_r^{\M}\circ f_{r+1}+ f_r\circ \(\pi_r^{\L}\circ f_{r+1}\circ \tau_r\)\circ \pi_r^{\L}\circ f_{r+1}+ f_r\circ \pi_r^{\Nh}\circ f_{r+1}\circ \tau_r\circ \pi_r^{\L}\circ f_{r+1}.$$Apply Lemma \tref{L3.15} to see that the expression inside the parentheses is the identity map, and apply 
(\tref{T3.40}) to see that $$f_r\circ m_{r+1}= f_r\circ \pi_r^{\Nh}\circ f_{r+1}\circ (-1+\tau_r\circ \pi_r^{\L}\circ f_{r+1}).$$ Thus, we have 
$$\split m_r\circ m_{r+1}& = \pi_{r-1}^{\M}(1-f_r\circ\tau_{r-1}\circ \pi_{r-1}^{\L})\circ f_r \circ m_{r+1}\\& = \pi_{r-1}^{\M}(1-f_r\circ\tau_{r-1}\circ \pi_{r-1}^{\L})\circ f_r\circ \pi_r^{\Nh}\circ f_{r+1}\circ (-1+\tau_r\circ \pi_r^{\L}\circ f_{r+1})\\& = \pi_{r-1}^{\M}\circ f_r\circ\[  \pi_r^{\Nh} -  \tau_{r-1}\circ \pi_{r-1}^{\L}\circ f_r\circ \pi_r^{\Nh}\] \circ  f_{r+1}\circ (-1+\tau_r\circ \pi_r^{\L}\circ f_{r+1})\endsplit$$  Apply Lemma \tref{L3.15}, once again, to see that the expression inside the brackets is zero; and therefore, (b) is established.

The interesting part of the proof of $\N\,_r\subseteq \ker\psi_r$ is
$$\split \psi_r\circ f_{r+1}\circ \pi_{r+1}^{\Nh} = { }&\(\fakeht \pi_r^{\M}-\pi_r^{\M}\circ f_{r+1}\circ \tau_r\circ \pi_r^{\L}\)\circ f_{r+1}\circ \pi_{r+1}^{\Nh} \\
{ } = { } &\pi_r^{\M}\circ f_{r+1}\circ\(\fakeht 1-\tau_r\circ \pi_r^{\L}\circ f_{r+1}\)\circ \pi_{r+1}^{\Nh}.\endsplit$$ Lemma \tref{L3.15} ensures that the expression inside the parentheses is zero.  
To prove that $\psi\:\F\to \M$ is a map of complexes, we must show that $$m_r\circ \psi_r(x_r)=\psi_{r-1}\circ f_r(x_r)\tag\tnum{37t12}$$ for all $x_r\in \F_r$. If $x_r\in \M\,_r$,  then the left side of (\tref{37t12}) is $$m_r(x_r)=\psi_{r-1}\circ f_r(x_r).$$ If $x_r\in \N\,_r$, then both sides of  (\tref{37t12}) are zero; and therefore, (\tref{37t12})  is established by (\tref{T3.18}). We complete the proof of (c) by identifying the kernel of $\psi$. 
  Let $x_r\in\Ker \psi_r$. Use (\tref{T3.18}) to write $x_r=y_r+z_r$ for some $y_r\in \M\,_r$ and some $z_r\in \N\,_r$. Observe that
$$0=\psi_r(x_r)=\psi_r(y_r)+\psi_r(z_r)=y_r.$$   We conclude that $\ker\psi_r=\N\,_r$ and the proof of (c) is complete.

We conclude by proving (d). It is clear that $\psi_r\circ \rho_r=\id$. The proof that $\rho$ is a map of complexes is much like the proof of (b). We see that $\rho_{r-1}\circ m_r$ is equal to 
$$\rho_{r-1}\circ \pi_{r-1}^{\M}\circ f_r - \rho_{r-1}\circ \pi_{r-1}^{\M}\circ f_r\circ \tau_{r-1}\circ \pi_{r-1}^{\L}\circ f_r.\tag\tnum{3str}$$The second term of (\tref{3str}) is equal to 
$$- \pi_{r-1}^{\M}\circ f_r\circ \tau_{r-1}\circ \pi_{r-1}^{\L}\circ f_r + 
\tau_{r-2} \circ \pi_{r-2}^{\L}\circ f_{r-1}\circ \pi_{r-1}^{\M}\circ f_r\circ \tau_{r-1}\circ \pi_{r-1}^{\L}\circ f_r.$$ Use (\tref{T3.40}) to see that the second term of (\tref{3str}) is
$$\split &- \pi_{r-1}^{\M}\circ f_r\circ 
\tau_{r-1}\circ \pi_{r-1}^{\L}\circ f_r 
-\[\tau_{r-2} \circ \pi_{r-2}^{\L}\circ f_{r-1}\fakeht \]\circ \pi_{r-1}^{\Nh}\circ f_r\circ \tau_{r-1}\circ \pi_{r-1}^{\L}\circ f_r\\ & -\tau_{r-2} \circ \pi_{r-2}^{\L}\circ f_{r-1}\circ \[\fakeht \pi_{r-1}^{\L}\circ f_r\circ \tau_{r-1}\]\circ \pi_{r-1}^{\L}\circ f_r.\endsplit$$ Lemma \tref{L3.15} guarantees that each bracketed expression is equal to the identity map. Apply (\tref{T3.40}) two times to see that the second term of (\tref{3str}) is equal to 
$$\split & \[\fakeht \pi_{r-1}^{\L}\circ f_r\circ \tau_{r-1}\]\circ \pi_{r-1}^{\L}\circ f_r - f_r\circ \tau_{r-1}\circ \pi_{r-1}^{\L}\circ f_r  +\[\tau_{r-2} \circ \pi_{r-2}^{\L}\circ f_{r-1}\fakeht\]\circ  \pi_{r-1}^{\Nh}\circ f_r\\& +\tau_{r-2} \circ \pi_{r-2}^{\L}\circ f_{r-1}\circ  \pi_{r-1}^{\M}\circ f_r.\endsplit $$ Once more, Lemma \tref{L3.15} ensures that each bracketed expression is the identity map. Thus,
the second term of (\tref{3str}) is 
$$ (\pi_{r-1}^{\L}+\pi_{r-1}^{\Nh}) \circ f_r - f_r\circ \tau_{r-1}\circ \pi_{r-1}^{\L}\circ f_r   +\tau_{r-2} \circ \pi_{r-2}^{\L}\circ f_{r-1}\circ  \pi_{r-1}^{\M}\circ f_r.$$ The definition of $\rho$ yields that
the first term of (\tref{3str}) is 
$$\pi_{r-1}^{\M}\circ f_r - \tau_{r-2} \circ \pi_{r-2}^{\L}\circ f_{r-1}\circ  \pi_{r-1}^{\M}\circ f_r.$$Combine the two most recent expressions to see that 
$$\rho_{r-1}\circ m_r= f_r-f_r\circ \tau_{r-1}\circ \pi_{r-1}^{\L}\circ f_r = f_r\circ ( \incl - \tau_{r-1}\circ \pi_{r-1}^{\L}\circ f_r)=f_r\circ \rho_r,$$ and the proof is complete. 
 $\qed$\enddemo

\demo{Proof of Theorem \tref{T3.9}} Assertions (a) and (c) are contained in 
Proposition \tref{P3.16}. To prove (b), we use Lemma \tref{L3.10} and the notation of Convention \tref{conv} to see that 
$$\allowdisplaybreaks \split m_1(1)^{(1)}\(\e_1^{[j]}\t \e_1^{[i]}\) & { } =
-u(\e_1^{[j]})\cdot v(\e_1^{[i]})+\[ \e_1^{[j]} \w (\W^{n-1}X) (\e_1^{[i]}[\eta])\](\eta)\\ & { } =-u_jv_i+(-1)^{\frac{n(n-1)}{2}}(\Adj \X)_{ij}\\ & { } = - (-1)^{\frac{n(n-1)}{2}}\[ (-1)^{\frac{n(n-1)}{2}} \v\u-\Adj \X\]_{ij}\\ m_1(2)^{(1)}(e_1^{[j]} \t 1)& { } = u[X(e_1^{[j]})] = \sum_iu_ix_{ij}\\ m_1(3)^{(1)}(e_1^{[i]} \t 1)& { } = [X(v)](e_1^{[i]}) = \sum_j x_{ij}v_j.\endsplit$$ A straightforward calculation using Definitions \tref{36D1} and \tref{42D1} shows that $m_0(\M\,_0)=0$, $$f_r(\M\,_r)\subseteq \(\fakeht I_1(\X)+I_1(\u)+I_1(\v)\)\cdot \F_{r-1},\quad\text{provided $r\neq 0$ and $r\neq 2$, and}$$ $f_2(\widehat{\M}\,_2)+m_2\(\[\F_2(3)^{(1)}\]''\)\subseteq \(\fakeht I_1(\X)+I_1(\u)+I_1(\v)\)\cdot \F_{1}$.  Assertion (d) is established and the proof is complete. 
$\qed$ \enddemo

The next result, which  used in the proof of Theorem \tref{T5.1}, is   a small piece of the fact that the complex $\M$ is self dual.  

\proclaim{Observation \tnum{R3.20}}If $(\M,m)$ is the complex of Theorem \tref{T3.9}, then $I_1(m_1)=I_1(m_{2n})$. \endproclaim

\demo{Proof}
The back of $\M$ looks like 
$$\eightpoint 0@>>>  \M\,_{2n}=\F_{2n}(3)^{(n)}= \W^nF\t \W^nF^* @>m_{2n}>>
\M\,_{2n-1}= \left\{ \matrix \[\F_{2n-1}(2)^{(n-1)}\]'= \W^{n-1}F\t \W^nF^* \\ \p \\ \[\F_{2n-1}(3)^{(n-1)}\]'=\W^{n-1}F\t \W^nF^*\\ \p\\\F_{2n-1}(4)^{(n-1)}=\W^{n-1}F^*\t \W^{n-1}F^*. \endmatrix \right.$$
Adopt the notation of Convention \tref{conv}.
One can readily check that 
$$m_{2n}(e_n\t \e_n)=\left\{ \matrix [X^*(u)](e_n)\t \e_n \in \[\F_{2n-1}(2)^{(n-1)}\]'\\+\\
[X(v)](e_n)\t \e_n \in \[\F_{2n-1}(3)^{(n-1)}\]'\\+\\
y \in\F_{2n-1}(4)^{(n-1)}, \endmatrix \right.$$
where
$$y= u(\e_n)\t v(\e_n)- \sum\limits_k (-1)^{k+1} e_1^{[k]}(\e_n)\t (\W^{n-1}X^*)\(e_{1}^{[1]}\w \dots\w \widehat{e_{1}^{[k]}}\w \dots \w e_{1}^{[n]}\).$$
For fixed integers $i$ and $j$, consider the homomorphism $$\f_{ij}\:\W^{n-1}F^*\t \W^{n-1}F^*\to \W^{n}F^*\t \W^{n}F^*,$$ which is given by
$$\f_{ij}(\a_{n-1}\t\b_{n-1})= \e_1^{[i]}\w\a_{n-1}\t \e_1^{[j]}\w\b_{n-1}.$$
A 
  short calculation, using Proposition \tref{A3}, yields  
$$ \eightpoint \f_{ij}(y)
=(-1)^{\frac{n(n-1)}{2}} \[(-1)^{\frac{n(n-1)}{2}}\v\u-\Adj \X \]_{ji} \cdot \e_n\t\e_n \in \W^nF^*\t\W^nF^*. \qed$$ \enddemo

\proclaim{Proposition \tnum{R3.21}} Adopt Data \tref{SU} with $n=2$, and let $\u,\X,\v$ be the matrices of Convention \tref{conv}.  If $\F$ is the complex   of Definition \tref{36D1}, then there is a split exact subcomplex $\N$ of $\F$ such that $\F/\N$ is the Koszul complex on the entries of $\Adj \X+\v\u$.\endproclaim

\demo{Proof}The proof is very similar to the proof of Theorem \tref{T3.9}. The main difference is due to the fact that $1=n-1$; and therefore, we must replace the modules of Definition \tref{D3.6} with
$$\eightpoint \[\F_2(2)^{(1)}\]'=\ker\mu_1,\quad \[\F_2(2)^{(1)}\]''=\im\ell_0,\quad
\[\F_2(3)^{(1)}\]'=\ker \lambda_1,\quad\text{and}\quad \[\F_2(3)^{(1)}\]''=\im \sigma_2.$$ Decompose $\F$ as $\L\p\M\p\Nh$, where $$\eightpoint \L=\left\{\phantom{\F_4(2)^{(2)}\ \p \ }\F_4(2)^{(2)}\ \p \ \matrix \F_3(2)^{(2)}\\\p\\ \F_3(3)^{(2)}\endmatrix \p \matrix \F_2(2)^{(2)}   \\\p\\ 
\F_2(3)^{(2)}\\\p\\ \[\F_2(3)^{(1)}\]''\\\p\\\F_2(4)^{(1)}\\\p\\\F_2(4)^{(0)}\endmatrix 
\ \p \ \matrix \F_1(1)^{(2)}\\\p\\ \F_1(1)^{(0)}\\\p\\ \F_1(2)^{(1)}\\\p\\\F_1(3)^{(1)}\\\p\\ \F_1(4)^{(0)}\endmatrix \ \p \ \matrix \F_0(1)^{(1)}\\\p\\ \F_0(1)^{(0)}\endmatrix\ \p  \ \F_{-1}(1)^{(0)}\right\},$$
$$\eightpoint \M=\left\{\phantom{\F_4(2)^{(2)}\ \p \ } \F_4(3)^{(2)}\ \p \   \F_3(4)^{(1)} \ \p\ \matrix \[\F_2(2)^{(1)}\]' \\\p\\ \[\F_2(3)^{(1)}\]' \endmatrix   \ \p \  \F_1(1)^{(1)}
  \ \p   \ 
\F_0(3)^{(0)} \phantom{\F_4(2)^{(2)}\ \p \ }\right\},$$and 
$$\eightpoint \Nh=\left\{ \F_5(4)^{(2)}\ \p \ \matrix \F_4(4)^{(2)}\\\p\\ \F_4(4)^{(1)}\endmatrix
\ \p\ \matrix \F_3(1)^{(2)}\\\p\\ \F_3(2)^{(1)}\\\p\\ \F_3(3)^{(1)}\\\p\\ 
\F_3(4)^{(2)}\\\p\\ \F_3(4)^{(0)}\endmatrix \ \p \ \matrix \F_2(1)^{(2)}\\\p\\
\F_2(1)^{(1)}\\\p\\ \[\F_2(2)^{(1)}\]''\\\p\\\F_2(2)^{(0)}\\\p\\\F_2(3)^{(0)}
\endmatrix \ \p \ \matrix \F_1(2)^{(0)}\\\p\\ \F_1(3)^{(0)}\endmatrix \ \p\ 
\F_0(2)^{(0)}\phantom{\F_4(2)^{(2)}\ \p \ }\right\}.$$Let $(\N,n)$ be the subcomplex of $\F$ which is given by 
$$\N\,_r=\Nh\,_r+f_{r+1}(\Nh\,_{r+1})\quad\text{and}\quad
 n_r=f_{r}|_{\N\,_r}.$$ 
In order to simplify the rest of the argument, we take $\ell_0(1)$ to be the element $e_1^{[1]}\t \e_1^{[1]}$ of $\F_2(2)^{(1)}$, and $\lambda_1\:\W^1F\t \W^1F^*\to \W^2F^*$ to be the map
$$\lambda_1\(\fakeht r_1\cdot e_1^{[1]}\t \e_1^{[1]} +
r_2\cdot e_1^{[1]}\t \e_1^{[2]} +
r_3\cdot e_1^{[2]}\t \e_1^{[1]} +
r_4\cdot e_1^{[2]}\t \e_1^{[2]} 
 \) =-r_1\cdot \e_2.$$ 
Define $\tau_r\:\L\,_r\to \Nh\,_{r+1}$ by $$\eightpoint\spreadmatrixlines{5pt}
\matrix\format\l\\  
\tau_{-1}(1)^{(0)}(1\t 1)=1\t 1\in\F_0(2)^{(0)}\\
\tau_{0}(1)^{(0)}(1\t \b_1)=- 1\t \b_1\in\F_1(2)^{(0)}\\
\tau_0(1)^{(1)}(\b_1\t 1)= 1\t \b_1\in \F_1(3)^{(0)}\\
\tau_{1}(1)^{(0)}(1\t \a_2)= 1\t \a_2\in\F_2(2)^{(0)}\\
\tau_{1}(1)^{(2)}(\a_2\t 1)=-1\t \a_2 \in\F_2(3)^{(0)}\\
\tau_1(2)^{(1)}(a_1\t 1)=\left\{\matrix \e_2\t a_1(\e_2)\in \F_2(1)^{(2)}\\+\\[v\w a_1](\e_2)\t \e_2\in \F_2(3)^{(0)}\endmatrix \right.\\
\tau_1(3)^{(1)}(a_1\t 1)=\left\{\matrix  a_1(\e_2)\t \e_2\in \F_2(1)^{(1)}\\+\\[ a_1\w u](\e_2)\t \e_2\in \F_2(2)^{(0)}\endmatrix \right.\\
\tau_1(4)^{(0)}(1\t 1)=\left\{\matrix  v_1\cdot \e_2\t 
 \e_1^{[2]} \in \F_2(1)^{(2)}\\+\\ v_1v_2\cdot 1\t \e_2\in\F_2(3)^{(0)}\\+\\ \ell_0(1) \in \[\F_2(2)^{(1)}\]''\endmatrix \right.\\
\endmatrix \hskip-6.5pt \matrix\format\l\\ 
\tau_2(2)^{(2)}(e_2\t 1)=\left\{\matrix  -\e_2 \t 1 \in \F_3(4)^{(2)}\\+\\ u\t\e_2 \in \F_3(3)^{(1)}\\+\\-\lambda_1[u\t u(\e_2)]\t\e_2\in\F_3(1)^{(2)}\endmatrix \right.\\
\tau_{2}(3)^{(1)}[\sigma_2(\e_2)]=\e_2 \t \e_2\in\F_3(1)^{(2)}\\\vspace{5pt}
\tau_2(3)^{(2)}(e_2\t 1)=\left\{\matrix  -1\t \e_2   \in \F_3(4)^{(0)}\\+\\ v\t\e_2 \in \F_3(2)^{(1)}\endmatrix \right.\\
\tau_{2}(4)^{(0)}(1\t \a_1)=\a_1(e_2) \t \e_2\in\F_3(2)^{(1)}\\
\tau_2(4)^{(1)}(\a_1\t 1)=\left\{\matrix  -\a_1(e_2)\t \e_2 \in \F_3(3)^{(1)}\\+\\ \lambda_1[\a_1(e_2)\t u(\e_2)]\t\e_2 \in \F_3(1)^{(2)}\endmatrix \right.\\
\tau_{3}(2)^{(2)}(e_2\t \a_1)=\e_2 \t \a_1\in\F_4(4)^{(2)}\\
\tau_{3}(3)^{(2)}(e_2\t \a_1)=\a_1\t \e_2 \in\F_4(4)^{(1)}\\
\tau_{4}(2)^{(2)}(e_2\t \e_2)=- \e_2\t \e_2 \in\F_5(4)^{(2)}.
\endmatrix
$$Use Definition \tref{42D1'}\,(b) and (c)  to define $\psi_r\:\F_r\to \M\,_r$ and $m_r\:\M\,_r\to \M\,_{r-1}$. It is not difficult to verify Lemma \tref{L3.15}. Proposition \tref{P3.16} is a formal result; and therefore, it also holds. A direct calculation now shows that   $\M$ is the Koszul complex on entries of $\Adj \X+\v\u$. Indeed, if we let 
$$\allowdisplaybreaks\eightpoint \gather g_1=x_{22}+v_1u_1,\ g_2=-x_{21}+v_2u_1,\ g_3=-x_{12}+v_1u_2,\ g_4=x_{11}+v_2u_2\ \text{in $R$,}\\
 w_1 =-\e_1^{[1]}\t \e_1^{[1]},\   
 w_2 =-\e_1^{[1]}\t \e_1^{[2]},\   
w_3 =-\e_1^{[2]}\t \e_1^{[1]},\   
 w_4 =-\e_1^{[2]}\t \e_1^{[2]},\ \text{in $\F_1(1)^{(1)}=\M\,_1$,}\\ 
 w_{12}=e_1^{[2]}\t \e_1^{[1]}\in \[\F_2(3)^{(1)}\]',\ 
 w_{13}=-e_1^{[2]}\t \e_1^{[1]}\in \[\F_2(2)^{(1)}\]',\ 
w_{14}=e_1^{[2]}\t \e_1^{[2]}\in \[\F_2(3)^{(1)}\]',\\ 
 w_{23}=\left\{\matrix e_1^{[1]}\t \e_1^{[1]}-e_1^{[2]}\t \e_1^{[2]} \in  \[\F_2(2)^{(1)}\]'\\ +\\
 -e_1^{[2]}\t \e_1^{[2]}\in \[\F_2(3)^{(1)}\]',\endmatrix\right.\\
 w_{24}=e_1^{[1]}\t \e_1^{[2]}\in \[\F_2(2)^{(1)}\]',\ 
 w_{34}=-e_1^{[1]}\t \e_1^{[2]}\in \[\F_2(3)^{(1)}\]'\ \text{in $\M\,_2$,}\\
w_{123}=\e_1^{[1]}\t \e_1^{[1]}, \ 
 w_{124}=\e_1^{[1]}\t \e_1^{[2]},\ 
 w_{134}=-\e_1^{[2]}\t \e_1^{[1]},\ 
 w_{234}=-\e_1^{[2]}\t \e_1^{[2]}\ \text{in $\F_3(4)^{(1)}=\M\,_3$, and}\\
w_{1234}=-e_2\t\e_2\in \F_4(3)^{(2)}=\M\,_4, 
\endgather$$then we see that
$$\eightpoint \gather m_1(w_i)=g_i,\ m_2(w_{ij})=g_i\cdot w_j-g_j\cdot w_i,\ m_3(w_{ijk})=g_i\cdot w_{jk} -g_j\cdot w_{ik}+ g_k\cdot w_{ij},\ \text{and}\\
m_4(w_{1234})= g_1\cdot w_{234} - g_2\cdot w_{134} + g_3\cdot w_{124} - g_4\cdot w_{123} \endgather
$$
for all $i$, $j$, and $k$. 
$\qed$\enddemo

\bigpagebreak

\SectionNumber=\xact\tNumber=1

\flushpar{\bf \number\SectionNumber.\quad Exactness.}

\medskip

\proclaim{Theorem \tnum{T4.1}} Fix an integer $n$, with $2\le n$.  Let $\u_{1\times n}$, $\X_{n\times n}$, and $\v_{n\times 1}$ be matrices of indeterminates over a commutative noetherian ring $R_0$,
$R$ be the polynomial ring $R_0[\{u_i, v_i, x_{ij}\mid 1\le i,j\le n\}]$, and  $u, X, v$ be the Data of \tref{SU} constructed from $\u,\X,\v$ by way of Convention \tref{conv}. If  $\F$ is the complex $\F\,[u,X, v]$ of Definition \tref{36D1},
  then the homology $H_i(\F)$ is zero for all integers $i$, except $i=0$.  \endproclaim

\demo{Proof}The proof proceeds by induction on $n$. If $n=2$, then the result is established in Proposition \tref{R3.21}. Henceforth, we assume that $3\le n$. The map $$f_0(2)^{(0)}\:\F_0(2)^{(0)}\to \F_{-1}(1)^{(0)}$$ is an isomorphism, and when this isomorphism is split from the complex $\F$, the resulting complex, $\overline{\F}$, has the same homology as $\F$ and looks like 
$$\overline{\F}\:\qquad 0\to\overline{\F}_{2n+1}\to \overline{\F}_{2n}\to\dots \to \overline{\F}_0\to 0.$$ Consequently, it suffices to apply  
 the acyclicity lemma \cite{\rref{BE75}, Corollary 4.2} and prove that the homology of the localization $\F_{x}$ is concentrated in position zero for each fixed indeterminate $x=x_{ij}$.   Let $R_0'$ be the ring  $R_0[\{x_{ik},x_{\ell j}\mid 1\le k\le n, 1\le \ell\le n, \ell\neq i\},\,x^{-1}]$.  It is easy to find matrices $M$ and $N$ with entries in $R_0'$ such that \roster
\item"{(a)}"  $\det M=\det N=1$,
\item"{(b)}"  $N\X M$ has the form $\[\smallmatrix 1&0\\0&\X'\endsmallmatrix\]$,
\item"{(c)}"  the entries of $X'$, $M^{-1}\v$, $\u N^{-1}$ form a sequence of indeterminates $T_1,\dots T_m$ over the ring $R_0'$, where $m=(n-1)^2+2n$, and 
\item"{(d)}"  the ring $R_{x}$ is equal to the polynomial ring $R_0'[T_1,\dots T_m]$.
\endroster Lemma \tref{L4.2}  guarantees that
$\F_{x}$ is isomorphic to the complex created using the data
$\u N^{-1}$, $N\X M$, $M^{-1}\v$. Thus, $\F_{x}$ is
isomorphic to the complex $(\G,g)$ of Lemma \tref{L4.5}. In the notation of  Lemma \tref{L4.5}, $\G$
is ``almost'' the total complex of $$0\to \F\,'@>\bmatrix
v_1\\-u_1\endbmatrix>> \F\,'\p\F\,'@>\bmatrix -u_1& -v_1\endbmatrix>> \F\,'\to
0.\tag\tnum{T4.30}$$ Indeed, $\G$ and the total complex of (\tref{T4.30})
differ only because the map $h_{r-2}$ in $g_r$ is not zero. Nonetheless,
the induction hypothesis, applied to the generic data of $\F\,'$, guarantees that the homology of
$\F\,'$ is concentrated in degree zero and that  $$0\to H_0(\F\,')@>\bmatrix
v_1\\-u_1\endbmatrix_* >>H_0( \F\,'\p\F\,')@>\bmatrix -u_1& -v_1\endbmatrix_*>>
H_0(\F\,')$$ is an exact sequence. (Keep in mind that $u_1$ and $v_1$ are indeterminates over the polynomial ring obtained by adjoining the entries of of the matrices which represent $X'$, $u'$ and $v'$ to $R'_0$.) Lemma \tref{L4.4} now yields that
$H_k(\F_{x})=0$ for all $k\neq 0$. 
$\qed$
\enddemo

\proclaim{Lemma \tnum{L4.2}}Adopt Data \tref{SU}. If $\theta\:F\to F$ is  an isomorphism with $\det \theta=1$, then 
 the complexes $$\F\,[\theta^{-1}(u), \theta^*\circ X,v], \quad \F\,[u,X,v], \quad\text{ and}\quad \F\,[u, X\circ \theta, \theta^{-1}(v)]$$ are  all isomorphic. 
  \endproclaim

\demo{Proof}Let $\overline{\F}=\F\,[\theta^{-1}(u), \theta^*\circ X,v]$,  $\F=\F\,[u,X,v] $,  and    $\Ft=\F\,[u, X\circ \theta, \theta^{-1}(v)]$.   We define maps $\Theta\:\F\to\Ft$  and $\Phi\: \F\to \overline{\F}$  by defining module isomorphisms
$$\Theta_r(i)^{(t)}\:\F_r(i)^{(t)}\to \Ft_r(i)^{(t)}\quad\text{and}\quad 
\Phi_r(i)^{(t)}\:\F_r(i)^{(t)}\to \overline{\F}_r(i)^{(t)}
,$$ for all $i$, $r$, and $t$, as follows: \vphantom{\tnum{T4.3} \tnum{T4.5} \tnum{T4.4} }
$$\eightpoint \matrix \format \l&\ \ \ \l\\ \Theta_r(1)^{(t)}\(  \a_t\t \b_s\)= \a_t\t (\W^s\theta^*)(\b_s),& \Phi_r(1)^{(t)}\(  \a_t\t \b_s\)=(\W^t\theta^*)(\a_t)\t \b_s,\\\vspace{4pt} 
\Theta_r(2)^{(t)}\(  a_t\t \b_s\)= (\W^t\theta^{-1})(a_t)\t (\W^s\theta^*)(\b_s),& \Phi_r(2)^{(t)} \( a_t\t \b_s\)= a_t\t \b_s,\\ \vspace{4pt}
\Theta_r(3)^{(t)}\(  a_t\t \b_s\)= a_t\t  \b_s,& \Phi_r(3)^{(t)}\(a_t\t\b_s\)= (\W^t\theta^{-1})(a_t)\t (\W^s\theta^*)(\b_s) ,\\ \vspace{4pt}
\Theta_r(4)^{(t)}\(  \a_t\t \b_s\)= \a_t\t (\W^s\theta^*)(\b_s),\ \text{and}& \Phi_r(4)^{(t)}\(  \a_t\t \b_s\)= (\W^t\theta^*)(\a_t)\t \b_s.\endmatrix  $$
A direct calculation  shows that $\Theta$  and $\Phi$ are maps of complexes; the following identities are used:
 $$\eightpoint \alignat 1 
 (\W^{s-1}\theta^*)[v(\b_s)] &= [\theta^{-1}(v)][(\W^{s}\theta^*)(\b_s)], \tag{\tref{T4.3}}\\
 [\W^{n-s}(X\circ \theta) ]\(\fakeht [(\W^{s}\theta^*)(\b_s)][\eta]\)& =  (\W^{n-s}X)(\b_s[\eta]), 
\quad \text{and}\tag{\tref{T4.5}}\\
 (\W^{t-1}\theta^{-1})\(\[\b_1 \w (\W^{n-t}X^*)(\a_t[\eta])\](\eta)\)&= \[\theta^*(\b_1 )\w [\W^{n-t}(X\circ \theta) ^*](\a_t[\eta])\](\eta).
\tag{\tref{T4.4}}\\ 
\endalignat$$
Identity (\tref{T4.3}) is obvious. To prove (\tref{T4.5}), apply Observation \tref{OBB}\,(a) to see that the left side is equal to $$(\W^{n-s}X  )\( \b_s\[(\W^n\theta)[\eta]\fakeht \]\).$$ The hypothesis $\det\theta=1$ ensures that 
$(\W^n\theta)[\eta]=e_n$. The right side of (\tref{T4.4}) is equal to 
$\((\W^{n-t+1}\theta^*)\[ \b_1 \w (\W^{n-t}X ^*)(\a_t[\eta])\]\)(\eta)$. Apply the isomorphism $\W^{t-1}\theta $ to both sides of the proposed identity and use Observation \tref{OBB}\,(a), as well as the hypothesis $\det\theta=1$, to complete the proof of (\tref{T4.4}). 
$\qed$\enddemo

 \proclaim{Lemma \tnum{L4.4}} Let $(\T,t)@>\a>>(\M,m)@>\b>>(\B,b)$ be a complex of complexes   over the commutative noetherian ring $R$,  $\{h_i:\T_i\to \B\,_{i+1}\}$ be a family of maps  which satisfies
$$b_{i+1}\circ h_i+h_{i-1}\circ t_i=0\quad \text{for all $i$,}$$ and  $(\D,d)$ be the complex with $\D\,_i=\B\,_i\p\M\,_{i-1}\p\T_{i-2}$ and $$d_i=\bmatrix b_i&(-1)^{i-1}\b_{i-1}&h_{i-2}\\0& m_{i-1}& (-1)^i \a_{i-2}\\ 0&0&t_{i-2}\endbmatrix.$$ If $$H_i(\T)=H_i(\M)=H_i(\B)=0\quad\text{for all $i\neq 0$, and}$$
$$0\to H_0(\T)@>\a_*>> H_0(\M) @>\b_*>> H_0(\B)  $$ is an exact sequence of $R-$modules,  then $H_i(\D)= 0$ for all   $i\neq 0$.\endproclaim 

\demo{Proof}Let $(\C, c)$ be the mapping cone of the $\a:\T\to\M$. In other words, $\C\,_i=\M\,_i\p\T_{i-1}$ and $$c_i=\bmatrix m_i &(-1)^{i-1}\a_{i-1} \\0&t_{i-1}\endbmatrix.$$ The long exact sequence of homology, 
which is associated to a mapping cone, yields   the exact sequence 
$$0@>>>H_0(\T)@>\a_*>>H_0(\M)@>>> H_0(\C)@>>>0,$$
as well as  $H_i(\C)=0$ for   $i\neq 0$. 
 Observe that $\D$ is the mapping cone of $[\b\ h]\: \C\to \B$. (The map $\C\,_i\to \B\,_i$ is given by  $[\b_i\ \ (-1)^ih_{i-1}]$.) Let $ \D\,'$ be the mapping cone of $[\b\ \ 0]\:\C\to \B$. The long exact sequence of homology gives  exact sequences
  $$\align &0\to H_1( \D\,')@>>>\H_0(\C)@>[\b\ 0]_*>>H_0(\B)@>>> H_0( \D\,')@>>>0, \quad\text{and}\\
&0\to H_1(\D)@>>>\H_0(\C)@>[\b\ h]_*>>H_0(\B)@>>> H_0(\D)@>>>0,\endalign
$$ as well as 
$H_i(\D)=H_i( \D\,')=0$ for $i\neq 0$ or $1$.  The hypothesis ensures that $[\b\ 0]_*$ is an injection.   The proof is complete because $[\b\ 0]_*$ and $[\b\ h]_*$ are the same map from $\H_0(\C)$ to $\H_0(\B)$. Indeed, if $\[\smallmatrix  z_0\\ z_{-1} \endsmallmatrix\]$ is a cycle in $\C\,_0$, then $t_{-1}( z_{-1})=0$ and $ z_{-1}=t_0(y_{0})$ for some $y_0$ in $\T_0$. It follows that $$h_{-1}( z_{-1})=h_{-1}\circ t_0(y_0)=-b_{1}\circ h_0(y_0)\in \im b_1. \qed$$ \enddemo

\remark{Remark \tnum{R4.12}} Adopt Data \tref{SU}. For each pair of integers $(i,j)$, let $f_r[i,j]$ be the composition
$$\F_r(j)@>\incl>> \F_r@>f_r>> \F_{r-1} @>\proj>> \F_{r-1}(i).$$ Notice that 
the map $f_r\:\F_r\to \F_{r-1}$ is given by 
$$f_r=\bmatrix f_r[1,1]&f_r[1,2]&f_r[1,3]& 0\\ f_r[2,1]&f_r[2,2]&0&f_r[2,4]\\
f_r[3,1]&0&f_r[3,3]& f_r[3,4]\\ 0&f_r[4,2]&f_r[4,3]&f_r[4,4]\endbmatrix. $$
Define a new map $\ft_r\:\F_r\to \F_{r-1}$ by $$\ft_r=\bmatrix f_r[1,1]&f_r[1,2]&f_r[1,3]& 0\\ -f_r[2,1]&f_r[2,2]&0&-f_r[2,4]\\
-f_r[3,1]&0&f_r[3,3]& -f_r[3,4]\\ 0&f_r[4,2]&f_r[4,3]&f_r[4,4]\endbmatrix. $$
Each map  $f_{\bullet}[i,i]$ is a Koszul map; and therefore, $f_{r}[i,i]\circ f_{r+1}[i,i]=0$. A quick calculation now shows that 
$$\dots @>\ft_{r+1}>> \F_r @>\ft_r>> \F_{r-1} @>\ft_{r-1}>> \dots \tag\tnum{T4.10}$$
is a complex, which we denote by $\Ft[u,X,v]$.\endremark  
 
\proclaim{Lemma \tnum{L4.7}} If the data of  \tref{SU} is adopted, then the complex $\F\,[u,X,-v]$ of Definition \tref{36D1} is isomorphic to the complex $\Ft[u,X,v]$ of $($\tref{T4.10}$)$.  
\endproclaim

\demo{Proof}For integers $i$, $r$, and $t$ define the module  isomorphism
$$\theta_r(i)^{(t)}\: \F_r(i)^{t}\to \F_r(i)^{(t)}\quad\text{by}$$
 $$ \eightpoint \theta_r(1)^{(t)}=(-1)^{r+1-t}\id, \  \theta_r(2)^{(t)}=(-1)^{r-t}\id, \  \theta_r(3)^{(t)}=(-1)^{t}\id, \ 
\text{and} \  \theta_r(4)^{(t)}=(-1)^{r-t}\id.$$It is not difficult to see that $\theta\: \Ft[u,X,v]\to \F\,[u,X,-v]$ is a homomorphism of complexes. $\qed$\enddemo

Most of the statement of the next result consists of notation. The only hypothesis is labeled (\tref{T4.6}). 
\proclaim{Lemma \tnum{L4.5}}Adopt Data \tref{SU}. 
Let $Re\p F'$ be a decomposition of $F$ into a direct sum of two free summands, and let $R\e\p F^{\prime *}$ be the corresponding decomposition of $F^*$. Let $u_1$ and $v_1$ be elements of $R$ and $u'$ and $v'$ be elements of $F'$ with $u=u_1e+u'$ and  $v=v_1e+v'$. Let $e_{n-1}'$ be the orientation element of $\W^{n-1}F'$ which has the property that $e_n=e\w e_{n-1}'$. 
Suppose that  $X$ may be decomposed as 
$$X=\bmatrix 1&0\\0&X'\endbmatrix,\tag\tnum{T4.6}$$ for some map $X'\:F'\to F'$.
Let $(\F,f)$ be the complex $\F\,[u,X,v]$ and $(\F\,',f')$ be the complex $\F\,[u', X', (-1)^{n-1}v']$. Then  there are homomorphisms $h_r\:\F\,'_r\to \F\,'_{r+1}$ such that 
the complexes $(\F,f)$ and $(\G,g)$ are isomorphic, where
$$\eightpoint \G_r=\F\,'_r\p\F\,'_{r-1}\p\F\,'_{r-1}\p\F\,'_{r-2}\quad\text{and}\quad
g_r=\bmatrix f_r'&(-1)^{r}u_1&(-1)^{r}v_1 &h_{r-2}\\ 0& f_{r-1}'&0&(-1)^rv_1\\
0&0&f'_{r-1}&(-1)^{r+1}u_1\\
0&0&0&f'_{r-2}\endbmatrix.$$
 \endproclaim

\demo{Proof}We have chosen to let $(\F\,',f')$ represent the complex $\F\,[u', X', (-1)^{n-1}v']$ because this choice of notation leads to a clean statement of the result; however, all of our calculations are made using the complex $\F\,[u', X',  v']$, which we refer to as 
 $(\F\,',f'')$. The maneuver from $(\F\,',f'')$ to $(\F\,',f')$ uses Lemma \tref{L4.7}, and occurs at the end of the proof. 
 
For each integer $r$, let  $h_r'\:\F\,'_r\to \F\,'_{r+1}$ be the homomorphism which is given by
$$\matrix \format \l \\ h_r'(1)^{(t)}(\a_{t}'\t \b_{s}')=
 \left\{\matrix
(-1)^{n-1+r}\[ (\W^{n-1-t}X^{\prime*})(\a'_{t}[e_{n-1}'])\](e_{n-1}')\t\b'_{s}\in \F\,'_{r+1}(2)^{(t)} \\+\\
(-1)^{n+1+r}\[(\W^{n-1-s}X')(  \b'_{s}[e_{n-1}'])\](e_{n-1}')\t\a'_{t}\in \F\,'_{r+1}(3)^{(s)},
 \endmatrix\right.\\ 
h_r'(2)^{(t)}(a'_t\t\b'_s)= (-1)^{r-1} (\W^{t}X')(a'_t)\t \b'_s\in\F\,'_{r+1}(4)^{(t)},\\
h_r'(3)^{(t)}(a'_t\t\b'_s)= (-1)^r \b'_s\t (\W^{t}X^{\prime*})(a'_t) \in\F\,'_{r+1}(4)^{(s)},\ \text{and}\\ h_r'(4)=0.\endmatrix$$

 The direct sum decompositions of $F$ and $F^*$ give rise to   decompositions
$$\F_r(i )^{(t)}=\left\{\matrix \F\,'_r(i )^{(t)}\\\p\\ R\e_{\text{left}}\t\F\,'_{r-1}(i )^{(t-1)}\\\p\\ R\e_{\text{right}}\t\F\,'_{r-1}(i)^{(t)}\\\p\\  R\e\t R\e\t\F\,'_{r-2}(i )^{(t-1)}\endmatrix\right.\qquad\qquad \F_r(j)^{(t)}= \left\{\matrix \F\,'_r(j)^{(t)}\\\p\\ Re\t \F\,'_{r-1}(j)^{(t-1)}\\\p\\ 
R\e\t \F\,'_{r-1}(j)^{(t)}\\\p\\ Re\t R\e\t \F\,'_{r-2}(j)^{(t-1)},\endmatrix\right.$$
for $i=1,4$ and $j=2,3$. The notation is self-explanatory; for example, if 
$t+s=r+1$, then 
$$\F_r(1)^{(t)}=\W^tF^*\t \W^{s}F^*= \(\W^tF^{\prime *}\p  [R\e\t\W^{t-1}F^{\prime *} ] \) \t  \(\W^sF^{\prime *}\p  [R\e\t\W^{s-1}F^{\prime *} ] \), $$and we let $R\e_{\text{left}}\t\F\,'_{r-1}(1)^{(t-1)}$ represent 
the summand $[R\e\t\W^{t-1}F^{\prime *}]\t \W^{s}F^{\prime *}$ of $\F_r(1 )^{(t)}$.
 Now that $\F_r$ has been decomposed into 16 summands, we recombine the pieces. Let 
$$\A\,_r=\F\,'_r(1)\p \F\,'_r(2)\p\F\,'_r(3)\p\F\,'_r(4),$$
$$\B\,_r=\(R\e_{\text{left}}\t\F\,'_{r-1}(1)\) \p\(Re\t \F\,'_{r-1}(2)\)\p \(R\e\t \F\,'_{r-1}(3)\) \p\(R\e_{\text{left}}\t\F\,'_{r-1}(4)\),$$
$$\C\,_r=\(R\e_{\text{right}}\t\F\,'_{r-1}(1)\) \p\(R\e\t \F\,'_{r-1}(2)\)\p \(Re\t \F\,'_{r-1}(3)\) \p\(R\e_{\text{right}}\t\F\,'_{r-1}(4)\),$$ and
$$\eightpoint \D\,_r=\(R\e\t R\e\t\F\,'_{r-2}(1)\) \p\(Re\t R\e\t \F\,'_{r-2}(2)\)\p \(Re\t R\e\t \F\,'_{r-2}(3)\) \p\(R\e\t R\e\t\F_{r-2}(4)\).$$
Observe that $\F_r=\A\,_r\p\B\,_r\p\C\,_r\p\D\,_r$. 

We establish the isomorphism $(\F,f)\iso (\G,g)$ in two steps. First we create a complex $(\Fh, \fh)$ and establish an isomorphism from $(\Fh,\fh)$ to  $(\F,f)$. Then, we prove that $(\Fh, \fh)$ is isomorphic to $(\G,g)$. 

For each pair of integers $(i,j)$, let  $f_r''[i,j]\:\F\,'_r(j)\to \F\,'_r(i)$ be the map which is described in Remark \tref{R4.12}. Let $\tri_r\ :\F\,'_r\to \F\,'_{r-1}$ be the map which is given by   $$\tri_r= \bmatrix f''_r[1,1]&f''_r[1,2]&f''_r[1,3]& 0\\ (-1)^{n-1}f''_r[2,1]&f''_r[2,2]&0&(-1)^{n-1}f''_r[2,4]\\
(-1)^{n-1}f''_r[3,1]&0&f''_r[3,3]&(-1)^{n-1} f''_r[3,4]\\ 0&f''_r[4,2]&f''_r[4,3]&f''_r[4,4]\endbmatrix. $$We notice that 
$$\tri_r=\cases f_r'',&\text{if $n$ is odd, and }\\\widetilde{f''}_r,&\text{(in the sense of Remark \tref{R4.12}), if $n$ is even.}\endcases$$ In any event, Remark \tref{R4.12} shows that $(\F\,', \tri)$ is a complex. We now define the complex $(\Fh, \fh)$ by   $\Fh_r=\G_r$, and 
$$\fh_r=\bmatrix  \tri_r&(-1)^{r}u_1&(-1)^{r}v_1 &h_{r-2}'\\ 0& \tri_{r-1}&0&(-1)^rv_1\\
0&0&\tri_{r-1}&(-1)^{r+1}u_1\\
0&0&0&\tri_{r-2}\endbmatrix.$$ We find it convenient to give an additional name to the components of $$\Fh_r=\F\,'_r\p\F\,'_{r-1}\p\F\,'_{r-1}\p\F\,'_{r-2}.$$ Let $\Ah\,_r$, $\Bh\,_r$, $\Ch\,_r$, and $\Dh\,_r$ represent $\F\,'_r$, the first  $\F\,'_{r-1}$, the second $\F\,'_{r-1}$, and $\F\,'_{r-2}$, respectively. 
There are natural isomorphisms
$$ \Ah\,_r@>\id>> \A\,_r,\quad \Bh\,_r@>\nat>> \B\,_r,\quad \Ch\,_r@>\nat>> \C\,_r,\quad \text{and}\quad \Dh\,_r@>\nat>> \D\,_r;$$ for example, if $a_t'\t\b_s'\in \F\,'_{r-1}(3)^{(t)}\subseteq \Bh\,_r$, then $$\nat( a_t'\t\b_s')=a_t'\t\e\w\b_s'\in R\e\t \F\,'_{r-1}(3)^{(t)}\subseteq \B\,_r.$$ 
 Let $\rho_r\:\F\,'_r\to \F\,'_r$ be the isomorphism 
$$\rho_r=\bmatrix 1&0&0&0\\0&-1&0&0\\0&0&-1&0\\0&0&0&1\endbmatrix\:\left\{\matrix  \F\,'_r(1)\\\p\\ \F\,'_r(2)\\\p\\\F\,'_r(3)\\\p\\ \F\,'_r(4)\endmatrix \right. @>\phantom{XXX}>>\left\{\matrix  \F\,'_r(1)\\\p\\ \F\,'_r(2)\\\p\\\F\,'_r(3)\\\p\\ \F\,'_r(4).\endmatrix \right.
 $$
For each integer $r$, consider the module isomorphism $\f_r\:\Fh_r\to \F_r$, which is given by
$$ \matrix\format\l&\qquad\qquad\l\\  \Ah\,_r=\F\,'_r@>\id>>\A\,_r,  &\Bh\,_r=\F\,'_{r-1} @>\rho_{r-1}>> \F\,'_{r-1}@>\nat>>\B\,_r,\\ \Ch\,_r=\F\,'_{r-1} @> (-1)^r>> \F\,'_{r-1}@>\nat>>\C\,_r,\ 
 \text{and} &\Dh\,_r=\F\,'_{r-2} @>(-1)^{r}\cdot\rho_{r-2} >>\F\,'_{r-2} @>\nat>>\D\,_r.\endmatrix$$
A very long, but straightforward, calculation yields that $$\f\: (\Fh,\fh)\to(\F,f)\  \text{is a homomorphism    of complexes.}\tag\tnum{T4.15}$$ 

Assume, for the time being, that (\tref{T4.15}) is established.
If $n$ is odd, then we take $h_r=h_r'$. In this case, $(\Fh, \fh)$ is already equal to $(\G,g)$ and the proof is complete. If $n$ is even, then 
let $\theta\:(\F\,',\tri)\to (\F\,',f')$ be the isomorphism of Lemma \tref{L4.7} and let $h_r=\theta_{r+1}\circ h_r'\circ \theta_r^{-1}$. It is not difficult to see that the isomorphism 
$$\bmatrix \theta_r&0&0&0\\0&\theta_{r-1}&0&0\\0&0&\theta_{r-1}&0\\ 0&0&0&\theta_{r-2}\endbmatrix\: \Fh_r\to \G_r$$ induces an isomorphism of complexes from   $(\Fh,\fh)$ to $(\G,g)$. Once again, the proof is complete.

Now, we turn our attention to proving (\tref{T4.15}). The proof involves sixteen calculations. We record four of these calculations and suppress the remaining twelve. No new ideas are required for the suppressed calculations. Fix integers $r$, $s$, and $t$ with $t+s=r+1$. We begin with $y=\a'_t\t\b'_s\in \F\,'_r(1)^{(t)}\subseteq \Ah\,_r$. Let $\Delta(\a'_t)=\sum\limits_j\a_1^{\prime[j]}\t\a_{t-1}^{\prime[j]}$ and
$\Delta(\b'_s)=\sum\limits_i\b_{1}^{\prime[i]}\t\b_{s-1}^{\prime[i]}$.  We show that
$$ f_r\circ \f_r (y)=\f_{r-1}\circ \fh_r(y).\tag\tnum{T4.16}$$ The left side of (\tref{T4.16}) is equal to 
$$ f_r(1)^{(t)} (\a'_t\t\b'_s)=  \left\{\matrix \a_t'\t v(\b'_s)\in \F_{r-1}(1)^{(t)}\\+\\ 
(-1)^r u(\a'_t)\t \b'_s \in \F_{r-1}(1)^{(t-1)}\\ +\\\sum\limits_i\[\b_1^{\prime[i]}\w (\W^{n-t}X^*)(\a_t'[\eta])\](\eta)\t\b_{s-1}^{\prime[i]}\in \F_{r-1}(2)^{(t-1)}\\+\\
 \sum\limits_j \[\a_1^{\prime[j]}\w(\W^{n-s}X)(\b'_s[\eta])\](\eta)\t\a_{t-1}^{\prime[j]}\in \F_{r-1}(3)^{(s-1)}.
 \endmatrix \right.$$It is clear that $v(\b'_s)=v'(\b'_s)$ and $u(\a'_t)=u'(\a'_t)$. Use Proposition \tref{A3}\,(a) to see that 
$$\a_t'(\eta)=(-1)^t e\w \a_t'(e_{n-1}').\tag\tnum{T4.20}$$
It follows that $$\split \[\b_1^{\prime[i]}\w (\W^{n-t}X^*)(\a_t'[\eta])\](\eta) &= 
(-1)^{t}\[\b_1^{\prime[i]}\w X(e)\w (\W^{n-t-1}X^{\prime *})(\a_t'[e_{n-1}'])\](e_{n})\\& = 
(-1)^{n-1}\[\b_1^{\prime[i]}\w (\W^{n-t-1}X^{\prime *})(\a_t'[e_{n-1}'])\](e_{n-1}');\endsplit$$ and we see that
the left side of (\tref{T4.16}) is equal to 
$$ \left \{\matrix \a_t'\t v'(\b'_s)\in \F\,'_{r-1}(1)^{(t)}\\+\\ 
(-1)^r u'(\a'_t)\t \b'_s \in \F\,'_{r-1}(1)^{(t-1)}\\ +\\
(-1)^{n-1}\sum\limits_i\[\b_1^{\prime[i]}\w (\W^{n-t-1}X^{\prime *})(\a_t'[e_{n-1}'])\](e_{n-1}')\t\b_{s-1}^{\prime[i]}\in \F\,'_{r-1}(2)^{(t-1)}\\+\\
 (-1)^{n-1}\sum\limits_j \[\a_1^{\prime[j]}\w(\W^{n-s-1}X')(\b'_s[e_{n-1}'])\](e_{n-1}')\t\a_{t-1}^{\prime[j]}\in \F\,'_{r-1}(3)^{(s-1)},
 \endmatrix \right.$$ which is the same as the right side of (\tref{T4.16}).
In our second calculation, we take $y$ equal to $\a'_{t-1}\t\b'_s\in \F\,'_{r-1}(1)^{(t-1)}\subseteq \Bh\,_r$. Let   $\Delta(\a'_{t-1})=\sum\limits_j\a_1^{\prime[j]}\t\a_{t-2}^{\prime[j]}$, and
$\Delta(\b'_s)=\sum\limits_i\b_{1}^{\prime[i]}\t\b_{s-1}^{\prime[i]}$.
The left side of (\tref{T4.16}) is equal to   $$\eightpoint  f_r(1)^{(t)}\(\e\w \a'_{t-1}\t\b'_s\fakeht\)= \left\{\matrix 
\e\w\a'_{t-1}\t v(\b_{s}')\in  \F_{r-1}(1)^{(t)} \\+\\
(-1)^r u[\e\w \a'_{t-1}]\t \b_s' \in \F_{r-1}(1)^{(t-1)}\\ +\\
 \sum\limits_i\[\b_1^{\prime [i]}\w (\W^{n-t}X^*)([\e\w\a_{t-1}'][\eta])\](\eta)\t\b_{s-1}^{\prime [i]}\in \F_{r-1}(2)^{(t-1)}\\+\\   
\[\e\w(\W^{n-s}X)(\b_s'[\eta])\](\eta)\t\a_{t-1}^{\prime}\in \F_{r-1}(3)^{(s-1)}
\\+\\
 -\sum\limits_j \[\a_1^{\prime [j]}\w(\W^{n-s}X)(\b_s'[\eta])\](\eta)\t\e\w\a_{t-2}^{\prime[j]}\in \F_{r-1}(3)^{(s-1)}.
 \endmatrix \right.$$ Use (\tref{T4.20}) to see that $\e\w(\W^{n-s}X)(\b_s'[\eta])=0$. It follows that the left side of (\tref{T4.16}) is equal to 
$$\eightpoint = \left\{\matrix 
\e\w\a'_{t-1}\t v'(\b_{s}')\in R\e_{\text{left}}\t\F\,'_{r-2}(1)^{(t-1)}\subseteq\B\,_{r-1} \\+\\
(-1)^r u_1\cdot \a'_{t-1}\t \b_s' \in \F\,'_{r-1}(1)^{(t-1)}\subseteq\A\,_{r-1}\\ +\\
(-1)^{r+1}\e\w u'(\a'_{t-1})\t \b_s' \in R\e_{\text{left}}\t\F\,'_{r-2}(1)^{(t-2)}\subseteq\B\,_{r-1}
\\+\\
(-1)^ne\w \sum\limits_i\[\b_1^{\prime [i]}\w (\W^{n-t}X^{\prime*})(\a_{t-1}'[e_{n-1}'])\](e_{n-1}')\t\b_{s-1}^{\prime [i]}\in Re\t\F\,'_{r-2}(2)^{(t-2)}\subseteq\B\,_{r-1}\\+\\   
 (-1)^n\sum\limits_j \[\a_1^{\prime [j]}\w(\W^{n-s-1}X')(\b_s'[e_{n-1}'])\](e_{n-1}')\t\e\w \a_{t-2}^{\prime[j]}\in R\e\t\F\,'_{r-2}(3)^{(s-1)}\subseteq\B\,_{r-1}.
 \endmatrix \right.$$
On the other hand,  the right side of (\tref{T4.16}) is equal to 
$$\eightpoint \f_{r-1}\left\{\matrix (-1)^r u_1\cdot \a_{t-1}'\t \b_{s}'\in \F\,'_{r-1}(1)^{(t-1)}\subseteq \Ah\,_{r-1}\\+\\ \tri_{r-1}(1)( \a_{t-1}'\t \b_{s}')\in \Bh\,_{r-1}\endmatrix \right.$$$$\eightpoint = \f_{r-1}\left\{\matrix (-1)^r u_1\cdot \a_{t-1}'\t \b_{s}'\in \F\,'_{r-1}(1)^{(t-1)}\subseteq \Ah\,_{r-1}\\+\\
f_{r-1}''[1,1]( \a_{t-1}'\t \b_{s}')\in \F\,'_{r-2}(1)\subseteq \Bh\,_{r-1}\\+\\
(-1)^{n-1}f_{r-1}''[2,1]( \a_{t-1}'\t \b_{s}')\in \F\,'_{r-2}(2)\subseteq \Bh\,_{r-1}\\+\\
(-1)^{n-1}f_{r-1}''[3,1]( \a_{t-1}'\t \b_{s}')\in \F\,'_{r-2}(3)\subseteq \Bh\,_{r-1}\endmatrix\right.$$
$$\eightpoint = \left\{\matrix (-1)^r u_1\cdot \a_{t-1}'\t \b_{s}'\in \F\,'_{r-1}(1)^{(t-1)}\subseteq\A\,_{r-1} \\+\\
\nat\circ f_{r-1}''[1,1]( \a_{t-1}'\t \b_{s}')\in R\e_{\text{left}}\t \F\,'_{r-2}(1) \subseteq\B\,_{r-1}\\+\\
(-1)^{n}\nat\circ f_{r-1}''[2,1]( \a_{t-1}'\t \b_{s}')\in Re\t\F\,'_{r-2}(2)\subseteq\B\,_{r-1} \\+\\
(-1)^{n}\nat\circ f_{r-1}''[3,1]( \a_{t-1}'\t \b_{s}')\in R\e\t \F\,'_{r-2}(3)\subseteq\B\,_{r-1}; \endmatrix\right.$$thus, (\tref{T4.16}) holds in this case. 
In our third calculation, we take  $y=\a'_t\t  \b'_{s-1}$ in $\F\,'_{r-1}(1)^{(t)}\subseteq \Ch\,_r$. Let 
$\Delta(\a'_t)=\sum\limits_j\a_1^{\prime[j]}\t\a_{t-1}^{\prime[j]}$, and
$\Delta(\b'_{s-1})=\sum\limits_i\b_{1}^{\prime[i]}\t\b_{s-2}^{\prime[i]}$.
The left side of (\tref{T4.16}) is equal to $ f_r(1)^{(t)}\((-1)^r\a'_t\t\e\w \b'_{s-1}\fakeht\)$
$$ \eightpoint= (-1)^r
\left\{\matrix 
\a'_t\t v(\e\w \b'_{s-1})\in \F_{r-1}(1)^{(t)} \\+\\
(-1)^r u(\a'_t)\t \e\w \b'_{s-1} \in \F_{r-1}(1)^{(t-1)}\\ +\\\sum\limits_i -\[\b_1^{\prime [i]}\w (\W^{n-t}X^*)(\a'_t[\eta])\](\eta)\t\e \w \b_{s-2}^{\prime [i]}\in \F_{r-1}(2)^{(t-1)}\\+\\
 \sum\limits_j \[\a_1^{\prime [j]}\w(\W^{n-s}X)([\e\w \b'_{s-1}][\eta])\](\eta)\t\a_{t-1}^{\prime [j]}\in \F_{r-1}(3)^{(s-1)}
 \endmatrix \right.$$
$$\eightpoint = (-1)^r\left\{\matrix v_1\cdot\a'_t\t  \b'_{s-1}\in \F\,'_{r-1}(1)^{(t)} \subseteq \A\,_{r-1}\\+\\
-\a'_t\t \e\w v'( \b'_{s-1})\in R\e_{\text{right}}\t\F\,'_{r-2}(1)^{(t)}\subseteq \C\,_{r-1} \\+\\
(-1)^r u'(\a'_t)\t \e\w  \b'_{s-1} \in R\e_{\text{right}}\t\F\,'_{r-2}(1)^{(t-1)}\subseteq \C\,_{r-1}\\ +\\(-1)^n\sum\limits_i \[\b_1^{\prime [i]}\w (\W^{n-t-1}X^{\prime*})(\a'_t[e_{n-1}'])\](e_{n-1}')\t \e\w \b_{s-2}^{\prime [i]}\in R\e\t\F\,'_{r-2}(2)^{(t-1)}\subseteq \C\,_{r-1}\\+\\
(-1)^n e\w \sum\limits_j \[\a_1^{\prime [j]}\w(\W^{n-s}X')( \b'_{s-1}[e_{n-1}'])\](e_{n-1}')\t\a_{t-1}^{\prime[j]}\in Re\t\F\,'_{r-2}(3)^{(s-2)}\subseteq \C\,_{r-1}.
 \endmatrix \right.$$ 
The right side of (\tref{T4.16}) is equal to 
$$\eightpoint \f_{r-1}\left\{\matrix (-1)^r v_1\cdot \a_{t}'\t \b_{s-1}'\in \F\,'_{r-1}(1)^{(t)}\subseteq \Ah\,_{r-1}\\+\\
f_{r-1}''[1,1]( \a_{t}'\t \b_{s-1}')\in \F\,'_{r-2}(1)\subseteq \Ch\,_{r-1}\\+\\
(-1)^{n-1}f_{r-1}''[2,1]( \a_{t}'\t \b_{s-1}')\in \F\,'_{r-2}(2)\subseteq \Ch\,_{r-1}\\+\\
(-1)^{n-1}f_{r-1}''[3,1]( \a_{t}'\t \b_{s-1}')\in \F\,'_{r-2}(3)\subseteq \Ch\,_{r-1}\endmatrix\right.$$
$$\eightpoint = \left\{\matrix (-1)^r v_1\cdot \a_{t}'\t \b_{s-1}'\in \F\,'_{r-1}(1)^{(t)}\subseteq\A\,_{r-1} \\+\\
(-1)^{r-1}\nat\circ f_{r-1}''[1,1]( \a_{t}'\t \b_{s-1}')\in R\e_{\text{right}}\t \F\,'_{r-2}(1) \subseteq\C\,_{r-1}\\+\\
(-1)^{r+n} \nat\circ f_{r-1}''[2,1]( \a_{t}'\t \b_{s-1}')\in R\e\t\F\,'_{r-2}(2)\subseteq\C\,_{r-1} \\+\\
(-1)^{r-n} \nat\circ f_{r-1}''[3,1]( \a_{t}'\t \b_{s-1}')\in Re\t \F\,'_{r-2}(3)\subseteq\C\,_{r-1}; \endmatrix\right.$$thus, (\tref{T4.16}) holds in this case. 
 In our fourth calculation, we take 
 $y=\a'_{t-1}\t\b'_{s-1}$ in $\F\,'_{r-2}(1)^{(t-1)}\subseteq \Dh\,_r$. Let $\Delta(\a'_{t-1})=\sum\limits_j\a_1^{\prime[j]}\t\a_{t-2}^{\prime[j]}$  and  $\Delta(\b'_{s-1})=\sum\limits_i\b_{1}^{\prime[i]}\t\b_{s-2}^{\prime[i]}$.
The left side of (\tref{T4.16}) is equal to  $f_r(1)^{(t)}\((1)^r\e\w \a'_{t-1}\t \e\w \b'_{s-1}\fakeht\)$
$$\eightpoint = (1)^r\left\{\matrix 
\e\w \a'_{t-1}\t v(\e\w \b'_{s-1})\in \F_{r-1}(1)^{(t)} \\+\\
(-1)^r u(\e\w \a'_{t-1})\t \e\w \b'_{s-1} \in \F_{r-1}(1)^{(t-1)}\\ +\\\[\e\w (\W^{n-t}X^*)([\e\w \a'_{t-1}][\eta])\](\eta)\t\b'_{s-1}\in \F_{r-1}(2)^{(t-1)}
\\ +\\-\sum\limits_i\[\b_1^{\prime [i]}\w (\W^{n-t}X^*)([\e\w \a'_{t-1}][\eta])\](\eta)\t\e\w \b_{s-2}^{\prime [i]}\in \F_{r-1}(2)^{(t-1)}
\\+\\
\[\e\w(\W^{n-s}X)([\e\w \b'_{s-1}][\eta])\](\eta)\t\a'_{t-1}\in \F_{r-1}(3)^{(s-1)}
\\+\\
 -\sum\limits_j \[\a_1^{\prime[j]}\w(\W^{n-s}X)([\e\w \b'_{s-1}][\eta])\](\eta)\t\e\w\a_{t-2}^{\prime[j]}\in \F_{r-1}(3)^{(s-1)},
 \endmatrix \right.$$which is equal to $(-1)^r$ times 
$$\eightpoint \left\{\matrix 
v_1\cdot   \e\w\a'_{t-1}\t  \b'_{s-1}\in R\e_{\text{left}}\t \F\,'_{r-2}(1)^{(t-1)} \subseteq \B\,_{r-1}\\+\\
-\e\w\a'_{t-1}\t \e\w v'(\b'_{s-1})\in R\e\t R\e\t\F\,'_{r-3}(1)^{(t-1)}\subseteq \D\,_{r-1} 
\\+\\
(-1)^r u_1\cdot \a'_{t-1}\t \e\w  \b'_{s-1} \in R\e_{\text{right}}\t\F\,'_{r-2}(1)^{(t-1)}\subseteq \C\,_{r-1}
\\ +\\
(-1)^{r+1}\e\w  u'(\a'_{t-1})\t \e\w  \b'_{s-1} \in R\e\t R\e\t\F\,'_{r-3}(1)^{(t-2)}\subseteq \D\,_{r-1}\\+\\
(-1)^{n-1}\[ (\W^{n-t}X^{\prime*})(\a'_{t-1}[e_{n-1}'])\](e_{n-1}')\t\b'_{s-1}\in \F\,'_{r-1}(2)^{(t-1)}\subseteq \A\,_{r-1}
\\ +\\(-1)^{n+1}e\w \sum\limits_i\[\b_1^{\prime [i]}\w (\W^{n-t}X^{\prime*})(
\a'_{t-1}[e_{n-1}'])\](e_{n-1}')\t\ \e\w \b_{s-2}^{\prime [i]}\in Re\!\t\! R\e\!\t
\!\F\,'_{r-3}(2)^{(t-2)} \subseteq \D\,_{r-1}  
\\+\\  
(-1)^{n+1}\[(\W^{n-s}X')(  \b'_{s-1}[e_{n-1}'])\](e_{n-1}')\t\a'_{t-1}\in \F\,'_{r-1}(3)^{(s-1)}\subseteq \A\,_{r-1}
\\+\\
 (-1)^{n+1}e\w \sum\limits_j \[\a_1^{\prime[j]}\w(\W^{n-s}X')(
\b'_{s-1}[e_{n-1}'])\](e_{n-1}')\t \e\w \a_{t-2}^{\prime[j]}\in Re\!\t\!R\e\!\t \!\F\,'_{r-3}(3)^{(s-2)}\subseteq \D\,_{r-1}.
 \endmatrix \right.$$The right side of (\tref{T4.16}) is equal to 
$$\eightpoint \f_{r-1}\left\{\matrix 
h_{r-2}'(1)^{(t-1)}(\a_{t-1}'\t \b_{s-1}')\subseteq \F\,'_{r-1}\subseteq \Ah\,_{r-1}\\+\\
(-1)^r v_1\cdot  \a_{t-1}'\t \b_{s-1}'\in   \F\,'_{r-2}(1)^{(t-1)}\subseteq \Bh\,_{r-1}\\+\\
(-1)^{r+1} u_1\cdot  \a_{t-1}'\t \b_{s-1}'\in   \F\,'_{r-2}(1)^{(t-1)}\subseteq \Ch\,_{r-1}\\+\\
f_{r-2}''[1,1]( \a_{t-1}'\t \b_{s-1}')\in \F\,'_{r-3}(1)\subseteq \Dh\,_{r-1}\\+\\
(-1)^{n-1}f_{r-2}''[2,1]( \a_{t-1}'\t \b_{s-1}')\in \F\,'_{r-3}(2)\subseteq \Dh\,_{r-1}\\+\\
(-1)^{n-1}f_{r-2}''[3,1]( \a_{t}'\t \b_{s-1}')\in \F\,'_{r-3}(3)\subseteq \Dh\,_{r-1}\endmatrix\right.$$
$$\eightpoint = \left\{\matrix
h_{r-2}'(1)^{(t-1)}(\a_{t-1}'\t \b_{s-1}')\subseteq \F\,'_{r-1}\subseteq \A\,_{r-1}\\+\\
(-1)^r v_1\cdot \e\w\a_{t-1}'\t \b_{s-1}'\in R\e_{\text{left}}\t \F\,'_{r-2}(1)^{(t-1)}\subseteq\B\,_{r-1} \\+\\
u_1\cdot \a_{t-1}'\t \e\w\b_{s-1}'\in R\e_{\text{right}}\t \F\,'_{r-2}(1)^{(t-1)}\subseteq\C\,_{r-1} \\+\\
(-1)^{r-1}\nat\circ f_{r-2}''[1,1]( \a_{t-1}'\t \b_{s-1}')\in R\e \t R\e\t \F\,'_{r-3}(1) \subseteq\D\,_{r-1}\\+\\
 (-1)^{n-1+r}\nat\circ f_{r-2}''[2,1]( \a_{t-1}'\t \b_{s-1}')\in Re\t R\e\t\F\,'_{r-3}(2)\subseteq\D\,_{r-1} \\+\\
 (-1)^{n-1+r}\nat\circ f_{r-2}''[3,1]( \a_{t-1}'\t \b_{s-1}')\in Re\t R\e\t \F\,'_{r-3}(3)\subseteq\D\,_{r-1}; \endmatrix\right.$$thus, (\tref{T4.16}) holds in this case. 
  $\qed$\enddemo

\bigpagebreak

\SectionNumber=\misc\tNumber=1

\flushpar{\bf \number\SectionNumber.\quad Further applications and questions.}

\medskip

\definition{Data \tnum{Dat2}} Fix an integer $n$, with $3\le n$.  Let $\u_{1\times n}$, $\X_{n\times n}$, and $\v_{n\times 1}$ be matrices of indeterminates over a commutative noetherian ring $R_0$,    $H$ be the ideal $H(\u,\X,\v)$ of Definition \tref{D1.2}  
 in the polynomial ring $R=R_0[\{u_i, v_i, x_{ij}\mid 1\le i,j\le n\}]$, and  $\overline{R}$ be the quotient $R/H$.
\enddefinition

\proclaim{Theorem \tnum{T5.1}} Adopt Data \tref{Dat2}. The ideal $H$ of $R$ is a perfect Gorenstein ideal of grade $2n$.\endproclaim

\demo{Proof}Let  $u,X,v$ be the data of \tref{SU} which is obtained from $\u,\X,\v$ by way of Convention \tref{conv}, and let $\M$ be the complex $\M\,[u,X,(-1)^{\frac{n(n-1)}{2}}v]$. Theorems \tref{T4.1} and \tref{T3.9} show that $\M$ is a resolution of $\overline{R}$ of length $2n$.  It follows from \cite{\rref{BE73w}} and Observation \tref{R3.20} that
$$ \gather \grade H\le \pd_R \overline{R}\le 2n\le \grade I_1(m_{2n})=\grade H\quad\text{and}\\ \Ext_R^{2n}(\overline{R},R)=H_0(\M^{\,*})=R/I_1(m_{2n})=\overline{R}. \qed\endgather $$ \enddemo

\remark{Remark} There are at least two other ways to calculate $\grade H$: one can calculate the height of $H$ as in \cite{\rref{V91}} (see also Lemma \tref{L5.1}), or one can specialize $H$, along the lines of Example \tref{E5.7}, and then calculate its grade. \endremark

\proclaim{Lemma \tnum{L5.1}} Adopt Data \tref{Dat2}. Let $s$ be $u_i$ for some $i$, or $v_i$ for some $i$, or some $n-2$ minor of $X$. Then, there exist indeterminates $Y_1,\dots Y_{n^2}$ such that  $\overline{R}_s=R_0[Y_1,\dots Y_{n^2}]_s$.
\endproclaim

\demo{Proof}Begin with $s=u_1$. It is not difficult to show that 
$$H_s= \( I_1(\u\X)+(\{(\v\u-\Adj \X)_{i1} \mid 1\le i\le n \})\fakeht\)R_s;$$ see, 
for example,  \cite{\rref{V91}, Proposition 3.3.2}. 
It follows that
$$\overline{R}_s=R_0[ \{x_{ij}\mid 2\le i\le n,\ 1\le j\le n\},\  u_1,\dots, u_n ]_s.$$  Now, we let $s$ be the determinant of the submatrix of $\X$ which is obtained by deleting rows and columns $n$ and $n-1$. Observe that 
$$\eightpoint H_s= \( \{(\u\X)_i\mid 1\le i\le n-2\},\ \{(\X\v)_i\mid 1\le i\le n-2\},\ \{(\v\u-\Adj \X)_{ij}\mid n-1\le i,j \le n \}\fakeht\)R_s;$$and therefore,
$$\overline{R}_s=R_0[ u_{n-1}, u_n, v_{n-1}, v_n,\ \{x_{ij}\mid (i,j)\neq (n-1,n-1), (n-1,n), (n,n-1)\ \text{or}\ (n,n) \} ]_s.   \qed$$\enddemo

\proclaim{Corollary \tnum{C5.2}} Adopt Data \tref{Dat2}. 
\roster 
\item"{(a)}" If $R_0$ is a domain, then so is  $\overline{R}$.
\item"{(b)}" Let $k$ be an integer with $k\le 8$.
\itemitem{$($i$)$} If $R_0$ satisfies the Serre condition $(S_{k+1})$, then so does $\overline{R}$.
\itemitem{$($ii$)$} If $R_0$ satisfies the Serre conditions $(R_k)$ and $(S_{k+1})$, then so does $\overline{R}$.
\endroster
In particular, if the ring $R_0$ is reduced, then so is $\overline{R}$; if the ring $R_0$ is normal, then so is $\overline{R}$. \endproclaim

\demo{Proof}Assertion (a) follows from  Theorem \tref{T5.1} together with \cite{\rref{V91}, Proposition 3.3.2}. A  version of this argument, which contains more details, may be found in \cite{\rref{BV88}, Theorem 2.10}. The proof of (b) also follows a standard argument; see, for example, \cite{\rref{KU90}, Theorem 9.4}. Let $P$ be a prime of $R$ with $H\subseteq P$ and $\depth \overline{R}_P\le k$. For (i) it suffices to show that $\overline{R}_P$ is Cohen-Macaulay; for (ii) it suffices to show that $\overline{R}_P$ is regular. Since $HR_P$ is a perfect ideal of grade $2n$ in the ring $R_P$, we know (use \cite{\rref{BV88}, Proposition 16.18}, if necessary) that
$$\grade PR_P-2n=\grade PR_P-\grade HR_P\le \grade\frac{PR_P}{HR_P}= \depth \overline{R}_P\le k.$$ It follows that $\grade P\le 2n +8$; thus, 
$I_1(\u)+I_1(\v)+I_{n-2}(\X)$ is not contained in $P$. It follows from Lemma \tref{L5.1} that $\overline{R}_P$ is a localization of a polynomial ring over the local ring $A=(R_0)_{P\cap R_0}$. The hypothesis $\depth \overline{R}_P\le k$ ensures that $\depth A\le k$. It follows that $A$ is Cohen-Macaulay in case (i) and regular in case (ii). The proof is complete.
$\qed$\enddemo

\proclaim{Corollary \tnum{C5.3}} Assume that the ring $R_0$ of  Data \tref{Dat2} is regular local with maximal ideal $\maxm$. If $\maxM$ is the maximal ideal $\maxm R+I_1(\u)+I_1(\v)+I_1(\X)$ of $R$,  then  the localization $\overline{R}_{\maxM}$ is not in the linkage class of a complete intersection.  
\endproclaim
\demo{Proof} Corollary \tref{C5.2} shows that $\overline{R}_{\maxM}$ satisfies the Serre condition $(R_8)$; however, Huneke and Ulrich \cite{\rref{HU87}, Theorem 4.2} have proved that a licci Gorenstein ring can satisfy $(R_7)$ only if it is a complete intersection. $\qed$ \enddemo 

\remark{Remark \tnum{R5.4}} A second proof of Corollary \tref{C5.3} is available when $R_0$ is a field and $n$ is odd.  In this case, the graded twists 
$$0\to \bigoplus\limits_j R(-d_{gj}) \to \dots \to \bigoplus\limits_j R(-d_{1j}) \to R$$ of the minimal resolution of $\overline{R}$  satisfy  the inequality 
$$ \max\limits_j\{d_{gj}\}\le (\grade H -1)\min\limits_j\{d_{1j}\}$$
of  \cite{\rref{HU87}, Cor. 5.13}, because Remark \tref{R3.1} shows that $\max\limits_j \{d_{gj}\}= n^2$, for $g=\grade H= 2n$, and if we take $\deg u=\deg v= \frac{n-1}{2}$, then  $\min\limits_j\{d_{1j}\}=\frac{n+1}{2}$. 
\endremark

The following consequence of Theorem  \tref{T5.1} is an application of the principle of the transfer of perfection; see, for example, \cite{\rref{BV88}, Theorem 3.5}.
 
\proclaim{Corollary \tnum{C5.6}} Let $\u'_{1\times n}$, $\X'_{n\times n}$, and $\v'_{n\times 1}$ be matrices with entries from a commutative noetherian ring $R'$, $H'$ be the ideal $H(\u',\X',\v')$ of Definition \tref{D1.2},
 $u', X', v'$ be the Data of \tref{SU} constructed from $\u',\X',\v'$ by way of Convention \tref{conv}, and   $\M$ be the complex $\M\,[u',X',(-1)^{\frac{n(n-1)}{2}}v']$. 
 If $H'$ is a proper ideal with $2n\le \grade H'$, then 
 $H'$ is a perfect Gorenstein ideal of grade equal to $2n$, and
 $\M$ is an $R'-$resolution of $R'/H'$. Furthermore, if $R'$ is a local ring with maximal ideal $\maxm$ and the entries of $\u'$, $\X'$, and $\v'$ are all in $\maxm$, then $\M$ is the minimal resolution of $R'/H'$.
 \endproclaim

\example{Example \tnum{E5.7}} Form  the ideal $H'=H(\u,\X',\v)$    using generic matrices $\u$ and $\v$ and a generic diagonal matrix $\X'=\diag(x_1,\dots, x_n)$. It is easy to see that $2n\le \grade H'$; and therefore Corollary \tref{C5.6} guarantees that $H'$ is a perfect Gorenstein ideal of grade $2n$. A significant amount of computer experimentation preceded the discovery of the resolutions $\F$ and $\M$. Among all   specializations of the generic Data
\tref{Dat2}, the ideals  of the present example, with $\deg u= [[\frac{n-1}{2}]]$, yield the best results when resolved using the computer program MACAULAY.
\endexample

As our final application, we determine which of the relations in the minimal resolution $\M$ are Koszul relations on the generators of the defining ideal of $H_0(\M)$. It is much easier to compute the subalgebra $k[\Tor_1]$ than it is to compute the entire algebra structure of $\Tor_{\bullet}$. Nonetheless, a significant amount of information is carried by this subalgebra; see, for example, \cite{\rref{Ku91}}.

\proclaim{Corollary \tnum{C5.10}}Adopt  the notation and hypotheses of Corollary \tref{C5.6}  with $3\le n$, $(R', \maxm, k)$ a local ring, and the entries of $\u'$, $\X'$, and $\v'$  all in $\maxm$. Then the subalgebra  $k[\Tor_1]$ of $\Tor_{\bullet}^{R'}(R'/H',k)$, which is generated by $\Tor_1$,   is isomorphic to the following quotient of the exterior algebra $\W_k^{\bullet}k^{n^2+2n}$ $:$ $$k[\Tor_1]\iso\frac{\W^{\bullet}_k(V_1\p V_2\p V_3)}{(\W^{n-1}V_2+ \W^{n-1}V_3+\W^2V_1+V_1V_2+V_1V_3+V_2V_3 )},$$ where $V_1$, $V_2$, and $V_3$ are vector spaces over $k$ with   $\dim V_1=n^2$ and $\dim V_2=\dim V_3=n$. Furthermore, there is an algebra embedding  
$$k[\Tor_1]\triv \(\Hom_k(k[\Tor_1],k) [-2n] \fakeht\) \hookrightarrow \Tor_{\bullet}^{R'}(R'/H',k).$$
\endproclaim

\demo{Proof} The final assertion reflects  the fact that $\Tor_{\bullet}^{R'}(R'/H',k)$ is a Poincar\'e algebra of length $2n$ (see \cite{\rref{BE77}, Theorem 1.5} or \cite{\rref{Av78}, Example 9.4}) and $k[\Tor_1]_{n}=0$. To complete the proof, we calculate $k[\Tor_1]$. Let $\overline{\phantom{xx}}$ represent the functor $\underline{\phantom{X}}\t_{R'}k$. If $\g\:\W^{\bullet}\M\,_1\to \F$ is any map of complexes which extends the commutative diagram
$$\CD \M\,_1 @>m_1>> \M\,_0 @>>> 0\\ @V \rho_1 VV @V \rho_0 VV @VVV\\ \F_1 @>f_1 >> \F_0 @> f_0 >> \F_{-1} 
\endCD\tag\tnum{D5.11}$$ of Proposition \tref{P3.16}\,(d),    and $\psi\:\F\to \M$ is the map of complexes from Definition \tref{42D1'}\,(b), then the composition $$\W^{\bullet}\M\,_1 @>\g>>\F@>\psi>>\M$$ induces an isomorphism
$$k[\Tor_1]\iso \frac{\W^{\bullet}\overline{\M\,_1}}{\Ker \overline {\psi}\circ \overline {\g}}.$$ Recall that $\M\,_1=\F_1(1)^{(1)}\p \F_1(2)^{(1)}\p\F_1(3)^{(1)}$. We define the relevant parts of $\g$ as follows:
\roster
\item"{(a)}"Let $r$ be an integer with $1\le r$, $a_1^{[1]}\t 1, \dots , a_1^{[r]}\t 1$ be elements from $\F_1(3)^{(1)}$, and  $a_r$ be the element $a_1^{[1]}\w\dots\w a_1^{[r]}$ of $\W^rF$. The map $\g_r\:\W^r\F_1(3)^{(1)}\to \F_r$ is given by 
$$ \eightpoint \g_r\(\fakeht (a_1^{[1]}\t 1)\w \dots \w( a_1^{[r]}\t 1)\)
=\left\{\matrix 1\t (\W^rX^*)(a_r)\in \F_r(2)^{(0)}
\\+\\ 
  a_r\t 1\in \F_r(3)^{(r)}.\endmatrix\right.$$
\item"{(b)}"Let  $r$ be an integer with $1\le r$, $a_1^{[1]}\t 1, \dots , a_1^{[r]}\t 1$ be elements from $\F_1(2)^{(1)}$,  and $a_r$ be the element $a_1^{[1]}\w\dots\w a_1^{[r]}$ of $\W^rF$. The map $\g_r\:\W^r\F_1(2)^{(1)}\to \F_r$ is given by 
$$ \eightpoint \g_r\(\fakeht (a_1^{[1]}\t 1)\w \dots \w( a_1^{[r]}\t 1)\)
=\left\{\matrix (-1)^{\frac{r(r-1)}{2}} a_r\t 1\in \F_r(2)^{(r)}\\+\\ 
(-1)^{\frac{r(r-1)}{2}}1\t (\W^rX)(a_r)\in \F_r(3)^{(0)}.\endmatrix\right.$$
 \item"{(c)}"If $\a_1\t\b_1\in \F_1(1)^{(1)}$ and $a_1\t 1\in \F_1(2)^{(1)}$, then
$$\g_2\(\fakeht [\a_1\t\b_1]\w [a_1\t 1]\) = \left\{\matrix -X(a_1)\w \a_1\t \b_1\in \F_2(1)^{(2)}\\+\\ -u(\a_1)\cdot a_1\t\b_1 \in \F_2(2)^{(1)}\\+\\ -v(\b_1)\cdot 1\t X(a_1)\w \a_1 \in \F_2(3)^{(0)} \\+\\ a_1(\b_1)\cdot \a_1\t 1 \in \F_2(4)^{(1)}. \endmatrix \right.$$
\item"{(d)}"If $\a_1\t\b_1\in \F_1(1)^{(1)}$ and $a_1\t 1\in \F_1(3)^{(1)}$, then
$$\g_2\(\fakeht [\a_1\t\b_1]\w [a_1\t 1]\) = \left\{\matrix
-\a_1 \t X^*(a_1)\w \b_1\in \F_2(1)^{(1)}\\+\\ u(\a_1)\cdot 1\t X^*(a_1)\w\b_1 \in \F_2(2)^{(0)}\\+\\ v(\b_1)\cdot a_1 \t \a_1 \in \F_2(3)^{(1)}\\+\\ a_1(\a_1)\cdot 1\t \b_1 \in \F_2(4)^{(0)}.  \endmatrix \right.$$
\item"{(e)}"If  $a_1\t 1\in \F_1(2)^{(1)}$ and $b_1\t 1\in \F_1(3)^{(1)}$, then
$$\g_2\(\fakeht [a_1\t 1]\w [b_1\t 1]\) = \left\{\matrix
-a_1 \t X^*(b_1)\in \F_2(2)^{(1)}\\+\\ -b_1\t X(a_1) \in \F_2(3)^{(1)}. 
 \endmatrix \right.$$
\item"{(f)}" The restriction of $\g_2$ to $\W^2\F_1(1)^{(1)}$ is the 
composition 
$$\W^2\F_1(1)^{(1)}@> s >> \F_1(1)^{(1)} \t \F_1(1)^{(1)} @> \widetilde{ \g_2} >> \F_2, $$ where $s$ is any splitting of the canonical map $\F_1(1)^{(1)}\t \F_1(1)^{(1)} \to \W^2\F_1(1)^{(1)}$, and  
$$\widetilde{ \g_2}\(\fakeht [\a_1\t\b_1]\t [\a_1'\t \b_1']\) = \left\{\matrix 
- v(\b_1')\cdot \a_1\w \a_1' \t \b_1 \in \F_2(1)^{(2)}\\+\\ -u(\a_1)\cdot \a_1'\t \b_1\w \b_1'\in \F_2(1)^{(1)}\\+\\  \[\b_1\w (\W^{n-2}X^*)[(\a_1\w\a_1')(e_n)]\](e_n)\t \b_1' \in \F_2(2)^{(1)}\\+\\
u(\a_1)u(\a_1')\t\b_1\w\b_1'\in\F_2(2)^{(0)}\\+\\
- \[\a_1'\w (\W^{n-2}X)[(\b_1\w\b_1')(e_n)]\](e_n)\t \a_1 \in \F_2(3)^{(1)}\\+\\
-v(\b_1)v(\b_1')\t\a_1\w\a_1'\in \F_2(3)^{(0)}.
 \endmatrix \right.$$
 \endroster It is not difficult to see that the part of $\g$ which we have defined may be extended to give a map of complexes $\g\: \W^{\bullet}\M\,_1\to \F$ which extends (\tref{D5.11}).
It is also easy to see that the kernel of $\overline {\psi}\circ \overline {\g}$ is equal to $$\eightpoint  \W^{n-1}\overline{\F_1(2)^{(1)}}+ \W^{n-1}\overline{\F_1(3)^{(1)}}+\W^2\overline{\F_1(1)^{(1)}}+\overline{\F_1(1)^{(1)}}\t \overline{\F_1(2)^{(1)}}+\overline{\F_1(1)^{(1)}}\t \overline{\F_1(3)^{(1)}}+\overline{\F_1(2)^{(1)}}\t \overline{\F_1(3)^{(1)}}, 
 $$ and the proof is complete. $\qed$ \enddemo

We conclude by recording some questions, which have not yet been addressed, about the ring $\overline{R}$ of Data \tref{Dat2}. In this discussion we   take $R_0=k$ to be the  field of complex numbers. Is $\overline{R}$ a Unique Factorization Domain? If not, what is its divisor class group? Is $\overline{R}$ rigid? What is the cotangent cohomology module $T^2(\overline{R}/k,\overline{R})=\Ext_{\overline{R}}^1(H/H^2,\overline{R})$? Does the minimal resolution of $\overline{R}$ admit the structure of a differential graded algebra? What is the rest of the algebra structure of
 $\Tor_{\bullet}^R(\overline{R},k)$? Is the Poincar\'e series
$$P_{\overline{R}}^k(z)=\sum_{i=0}^{\infty}\dim_k\Tor_{\overline{R}}(k,k) z^i$$ a rational function? Are there interesting ideals, which are analogous to $H$, when the matrix $\X$ is not square? 

\end
\Refs

\ref \no \rnum{Av78} \by L\. Avramov \paper Small homomorphisms of
local rings
\jour J\. Alg\. \yr 1978 \vol 50 \pages 400--453 \endref 

\ref\no\rnum{Bss} \by H. Bass \paper On the ubiquity of Gorenstein rings \jour Math. Z. \vol 82 \yr 1963 \pages 8--28\endref 

\ref \no \rnum{BV88} \by W\. Bruns and U\. Vetter \book Determinantal
rings
\bookinfo Lecture Notes in Mathematics {\bf 1327} \yr 1988 \publ Springer
Verlag \publaddr Berlin Heidelberg New York \endref 

\ref \no \rnum{BE73w} \by D\. Buchsbaum and D\. Eisenbud \paper What makes a
complex
exact? \jour 
J\.
Alg\. \vol 25 \yr 1973 \pages 259--268 \endref



\ref
\no \rnum{BE75}
\by D\. Buchsbaum and D\. Eisenbud
\paper Generic free resolutions and a family of generically perfect ideals
\jour Advances Math\.
\yr 1975
\vol 18
\pages 245--301
\endref

\ref \no \rnum{BE77} \by D\. Buchsbaum and D\. Eisenbud \paper Algebra
structures for finite free resolutions, and some structure theorems for ideals
of codimension 3 \jour Amer\. J\. Math\. \yr 1977 \vol 99 \pages 447--485
\endref  



\ref \no \rnum{HU87} \by C\. Huneke and B\. Ulrich \paper The structure
of
linkage \jour Annals of Math\. \yr 1987 \vol 126 \pages 277--334 \endref 

\ref\no \rnum{Ku91} \by A\. Kustin \paper Classification of the
Tor$-$algebras of codimension four almost complete intersections
\jour Trans. Amer. Math. Soc. \vol339 \yr 1993 \pages 61--85 \endref



\ref \no \rnum{KU90} \by A\. Kustin and B\. Ulrich \paper A family of
complexes
associated to an almost alternating map, with applications to residual
intersections \jour Mem\. Amer\. Math\. Soc\. \yr 1992
\vol 95 \pages 1--94 \endref 

\ref\no \rnum{V91}\by W\. Vasconcelos\paper On the equations of Rees
algebras
\jour J\. reine angew\. Math\. \vol 418 \yr 1991\pages 189--218 \endref 

\endRefs
